\documentclass{amsart}

\newcommand{\R}{\ensuremath{\mathbb{R}}}

\newcommand{\matr}[4]{
\left[\begin{array}{cc}
#1&#2\\
#3&#4
\end{array}
\right]}

   \newcommand{\tp}[3]{\left\{#1#2#3\right\}}

     \newcommand{\tpc}[3]{\left\{#1,#2,#3\right\}}

\theoremstyle{remark}

\newcommand{\CC}{{\bf C}}

\newtheorem{lemma}{Lemma}[section]
\newtheorem{theorem}[lemma]{Theorem}
\newcommand{\pf}{\noindent{\it Proof}.}
\newtheorem{proposition}[lemma]{Proposition}
\newtheorem{remark}[lemma]{Remark}
\newtheorem{definition}[lemma]{Definition}
\newtheorem{corollary}[lemma]{Corollary}
\newtheorem{problem}{Problem}

\newcommand{\kbta}{K\overline{\otimes}A}
\newcommand{\kta}{K\otimes A}

\begin{document}
\title{A holomorphic  characterization of  operator algebras}
\author{Matthew Neal}
\address{Department of Mathematics, Denison University, Granville,
Ohio 43023}
   \email{nealm@denison.edu}

\author{Bernard Russo}
\address{Department of Mathematics, University of California,
Irvine, California 92697-3875}
   \email{brusso@math.uci.edu}

\subjclass{Primary  47L25, 47L70  Secondary 17C65}

\date{}


\keywords{operator space, operator algebra, complete holomorphic vector field, partial triple product, symmetric part}

\begin{abstract}
A necessary and sufficient condition for an operator space to support a
 multiplication making it completely isometric and isomorphic  to a unital
 operator algebra is proved. The condition involves only the holomorphic
 structure of the Banach spaces underlying the operator space.
\end{abstract}

\maketitle

\section{Introduction and background}

\subsection{Introduction}
If $A$ is an operator algebra, that is, an associative subalgebra of $B(H)$, then $M_n(A)=A\otimes M_n(\CC)$ is a subalgebra of  $B(\oplus_1^n H)$ and its multiplication is contractive, that is, $\|XY\|\le \|X\|\|Y\|$ for $X,Y\in M_n(A)$, where $XY$ denotes the matrix or operator product of $X$ and $Y$. Conversely, if an operator space $A$ (i.e., a closed linear subspace of $B(H)$) is also a unital (not necessarily associative) Banach algebra  with respect to a product $x\cdot y$ which is completely contractive in the above sense,  then according to  \cite{BleRauSin90}, it is completely isometric via an algebraic isomorphism to an operator algebra (i.e., an associative subalgebra of some $B(K)$).   Our main result (Theorem~\ref{thm:main2}) drops the algebra assumption on $A$  in favor of a holomorphic assumption. Using only natural conditions on holomorphic vector fields on Banach spaces, we are able to construct an algebra product on $A$ which is completely contractive and unital, so that the result of \cite{BleRauSin90} can be applied. 
Thus we give a holomorphic
characterization of operator spaces which are completely isometric to operator
algebras. This paper is a companion to \cite{NeaRus03}  where the authors gave 
holomorphic characterizations of operator spaces that are completely isometric to
a C$^*$-algebra or to a ternary ring of operators (TRO).

This paper  is also an instance where the consideration of  a ternary product, called the {\it partial triple product}, which arises from the holomorphic structure via the so called {\it symmetric part}  of the Banach space (defined in subsection 1.2), leads to results for binary products.  Examples of this phenomenon occurred in  \cite{ArazyvN}, \cite{AraSol90} where this idea is used to describe the algebraic properties of  isometries of certain operator algebras.  The method was also used in 
\cite{KauUpm78} to show that Banach spaces with holomorphically equivalent unit balls are linearly isometric (see \cite{Arazysurvey} for an exposition of \cite{KauUpm78}).   Another example is \cite{Kaneda04}, where it is proved that all operator algebra products on an operator space $A$ are of the form $x \cdot y = xa^{\ast}y$ for an element $a$ which  lies in the injective envelope $I(A)$. Here the ``quasimultiplier" $a$ lies in the symmetric part of $I(A)$.

Our technique is to use a variety of elementary  isometries on $n$ by $n$ matrices over $A$ (most  of the time, $n=2$) and to exploit the fact that isometries of arbitrary Banach spaces preserve the partial triple product. The first occurrence  of this technique appears in section 2, where for each $n$ a contractive projection $P_n$ on $\kbta$ ($K$= compact operators on separable infinite dimensional Hilbert space) with range $M_n(A)$ is constructed  as a convex combination of isometries. We define the {\it completely symmetric part} of $A$ to be the intersection of $A$ (embedded in $\kbta$) with the symmetric part of $\kbta$ and show 
it is  the image under $P_1$ of the symmetric part of $\kbta$. It follows from \cite{NeaRus03} that the completely symmetric part of $A$ is completely isometric to a TRO, which is a crucial tool in our work. 

We note that  if $A$ is a subalgebra of  $B(H)$ containing the identity operator $I$, then by \cite[Cor. 2.9(i)]{AraSol90}, its symmetric part is the maximal C$^*$-subalgebra $A\cap A^*$ of $A$. For the same reason, the symmetric part of the operator algebra $\kbta$ is the maximal $C^*$-subalgebra of $K\overline{\otimes}\  B(H)$ contained in $\kbta$, namely $\kbta\cap(\kbta)^*$, which shows that the completely symmetric part of $A$ coincides with  the symmetric part $A\cap A^*$ of $A$, and therefore contains $I$.  Moreover, by  \cite[Cor. 2.9(ii)]{AraSol90}, the partial triple product in $M_n(A)$ is the restriction of the triple product on $M_n(B(H))$.  

Thus the conditions (i) and (ii)  in our main theorem, which is stated here as Theorem~\ref{thm:1.1}, hold when $A$ is an operator algebra. This is evident from the restatement of (i) and (ii)  in (\ref{eq:0706121}).
In this theorem, for any element $v$ in the symmetric part  of a Banach space $X$, $h_v$ denotes the corresponding complete holomorphic vector field on the open unit ball of $X$. (Complete holomorphic vector fields and the symmetric part of a Banach space are recalled in subsection 1.2.)
\medskip

\begin{theorem}\label{thm:1.1}
An operator space $A$ is completely isometric to a unital operator
algebra if and only there exists an element $v$  of norm one in the completely symmetric part of $A$ such that:
\begin{description}
\item[(i)] $h_v(x+v)-h_v(x)-h_v(v)+v=-2x$ for all $x\in A$
\smallskip
\item[(ii)]  Let  $V=\hbox{diag}(v,\ldots,v)\in M_n(A)$.  For all $X\in M_n(A)$,
$$
\|V-h_{V}(X)\|\le \|X\|^2.
$$
\end{description}
\end{theorem}
\smallskip

Although we have phrased this theorem in holomorphic terms, it should be noted that the two conditions can be restated in terms of partial triple products as
\begin{equation}\label{eq:0706121}\{xvv\}=x\hbox { and } 
\|\{XVX\}\|\le \|X\|^2.
\end{equation}

Let us consider another  example.
Suppose that $A$ is a TRO, that is, a closed subspace of $B(H)$ closed under the ternary product $ab^*c$. Since $K\overline{\otimes}B(H)$ is a TRO, hence a $JC^*$-triple, it is equal to its symmetric part, which shows that the completely symmetric part of $A$ coincides with $A$.

Now suppose that the TRO $A$ contains an element $v$ satisfying $xv^*v=vv^*x=x$ for all $x\in X$.  Then it is trivial that $A$ becomes a unital C$^*$-algebra for the product $xv^*y$, involution $vx^*v$, and unit $v$.  By comparison, our main result starts only with an operator space $A$ containing a distinguished element $v$ in its completely symmetric part  having a unit-like property. We then construct a binary product from a property of  the partial triple product induced by the holomorphic structure.  $A$, with this binary product,  is then shown to be completely isometric to a unital operator algebra. The first assumption is unavoidable since the result of \cite{BleRauSin90} fails in the absence of a unit element.

According to \cite{BleZarPNAS},  ``The one-sided multipliers of an operator space $X$ are a key to the `latent operator algebraic structure' in $X$.''  The unified approach through multiplier operator algebras developed in \cite{BleZarPNAS} leads to simplifications of known results and applications to quantum $M$-ideal theory.  They also state
``With the extra structure consisting of the additional matrix norms on an operator algebra, one might expect to not have to rely  as heavily on other structure, such as the product.''  Our result is certainly in the spirit of this statement. 
\smallskip

In the rest of this section, a review of  operator spaces, Jordan triples, and holomorphy is given.
The completely symmetric part of an arbritary operator space $A$ is defined in section 2.  The binary product $x\cdot y$  on $A$  is constructed in section 3 using  properties of isometries on 2 by 2 matrices over $A$ and it is shown that the symmetrized product can be expressed in terms of the partial Jordan triple product  on $A$. 
It is worth noting that only the first hypothesis in Theorem~\ref{thm:1.1} is needed to prove the existence and properties of the binary product $x\cdot y$.  \smallskip

\subsection{Background}

In this section, we recall some basic facts that we use about operator spaces, Jordan triples, and holomorphy in Banach spaces. Besides the sources referenced in this section, for more facts and details on the first two topics, see \cite{U1},\cite{U2},\cite{Chubook}  and \cite{EffRau00},\cite{Paulsen02},\cite{Pisier03},\cite{BleLeM04},  respectively.

By an {\bf operator space}, sometimes called a quantum Banach space, we mean a closed linear subspace $A$ of $B(H)$ for some complex Hilbert space $H$, equipped with the matrix norm structure obtained by the identification of $M_n(B(H))$
with $B(H\oplus H\oplus\cdots\oplus H)$.  Two operator spaces are {\bf completely isometric} if there is a linear isomorphism between them which, when applied elementwise to the corresponding spaces of $n$ by $n$ matrices, is an isometry for every $n\ge 1$. 

By an {\bf operator algebra}, sometimes called a quantum operator algebra, we mean a closed associative subalgebra $A$ of $B(H)$, together with its matrix norm structure as an operator space.

One important example of an operator space is a {\bf ternary ring of operators}, or TRO, which is an operator space in $B(H)$ which contains $ab^*c$ whenever it contains $a,b,c$.

A TRO is a special case of a {\bf $JC^*$-triple}, that is, a closed subspace of $B(H)$ which contains the symmetrized ternary product $ab^*c+cb^*a$ whenever it contains $a,b,c$.  More generally, a {\bf $JB^*$-triple} is a complex Banach space equipped with a triple product $\{x,y,z\}$ which is linear in the first and third variables, conjugate linear in the second variable, satisfies the algebraic identities
\[
\{x,y,z\}=\{z,y,x\}
\]
and
\begin{equation}\label{eq:main}
\{a,b,\{x,y,z\}\}=\{\{a,b,x\},y,z\}-\{x,\{b,a,y\},z\}+\{x,y,\{a,b,z\}\}
\end{equation}
and the analytic conditions that the linear map $y\mapsto\{x,x,y\}$ is hermitian and positive and $\|\{x,x,x\}\|=\|x\|^3$.

The following two theorems are needed in what follows.

\begin{theorem}[Kaup \cite{Kaup83}]\label{thm:kaup83}
The class of $JB^*$-triples coincides with the class of complex Banach spaces whose open unit ball is a bounded symmetric domain.
\end{theorem} 

\begin{theorem}[Friedman-Russo \cite{FriRus85}, Kaup \cite{Kaup84}, Stacho \cite{Stacho82}]\label{thm:proj}
The class of $JB^*$-triples is stable under contractive projections.  More precisely, if $P$ is a contractive projection on a $JB^*$-triple $E$ with triple product denoted by $\{x,y,z\}_E$, then $P(E)$ is a $JB^*$-triple with triple product given by $\{a,b,c\}_{P(E)}=P\{a,b,c\}_E$ for $a,b,c\in P(E)$.
\end{theorem}

For a $JB^*$-triple, the following identity is a consequence of the Gelfand Naimark Theorem (\cite[Corollary 3]{FriRus86}):
$$
\|\{xyz\}\|\le \|x\|\|y\|\|z\|.
$$
This suggests Problem~\ref{prob:1} at the end of this paper.

The following two theorems, already mentioned above,  are instrumental in this work.

\begin{theorem}[Blecher,Ruan,Sinclair \cite{BleRauSin90}]\label{thm:brs}
If an operator space supports a unital Banach algebra structure in which the product (not necessarily associative) is completely contractive, then the operator space is completely isometric to an operator algebra.
\end{theorem}

\begin{theorem}[Neal,Russo \cite{NeaRus03}]\label{thm:tro}
If an operator space has the property that the  open unit ball of the space of $n$ by $n$ matrices is a bounded symmetric domain for every $n\ge 2$, then the operator space is completely isometric to a TRO.
\end{theorem}

Finally, we review the construction and properties of the partial Jordan triple product in an arbitrary Banach space.  Let $X$ be a complex Banach space with open unit ball $X_0$. Every holomorphic function $h:X_0\rightarrow X$, also called a holomorphic vector field,  is locally integrable, that is, the initial value problem
\[
\frac{\partial}{\partial t}\varphi(t,z)=h(\varphi(t,z))\ ,\ \varphi(0,z)=z, 
\]
has a unique solution for every $z\in X_0$ for $t$ in a maximal open interval $J_z$ containing 0.  A {\bf complete holomorphic vector field} is one for which $J_z=\R$ for every $z\in X_0$.

It is a fact that every complete holomorphic vector field is the sum of the restriction of a skew-Hermitian bounded linear operator $A$ on $X$
and a function $h_a$ of the form
$
h_a(z)=a-Q_a(z),
$
where $Q_a$ is a quadratic homogeneous polynomial on $X$.  

The {\bf symmetric part} of $X$ is the orbit of 0 under the set of complete holomorphic vector fields, and is denoted by  $S(X)$.  It is a closed subspace of $X$ and is equal to $X$ precisely when $X$ has the structure of a $JB^*$-triple
(by Theorem~\ref{thm:kaup83}).

If $a\in S(X)$, we can obtain a symmetric bilinear form on $X$, also denoted by $Q_a$ via the polarization formula
\[
Q_a(x,y)=\frac{1}{2}\left(Q_a(x+y)-Q_a(x)-Q_a(y)\right)
\]
and then the partial Jordan triple product  $\{\cdot,\cdot,\cdot\}:X\times S(X)\times X\rightarrow X$ is defined by $\{x,a,z\}=Q_a(x,z)$.
The space $S(X)$ becomes a $JB^*$-triple in this triple product.

It is also true that the ``main identity'' (\ref{eq:main}) holds whenever $a,y,b\in S(X)$ and $x,z\in X$.
The following lemma is an elementary consequence of the definitions.

\begin{lemma}\label{lem:isometry}
If $\psi$ is a linear isometry of a Banach space $X$ onto itself, then
\begin{description}
\item[(a)] For every complete holomorphic vector field $h$ on $X_0$, $\psi\circ h\circ \psi^{-1}$ is a complete holomorphic vector field.  In particular, for $a\in S(X)$,
$\psi\circ h_a\circ \psi^{-1}=h_{\psi(a)}$.
\item[(b)] $\psi(S(X))=S(X)$ and $\psi$ preserves the partial Jordan triple product:
\[
\psi\{x,a,y\}=\{\psi(x),\psi(a),\psi(y)\}\hbox{ for }a\in S(X),\ x,y\in X.
\]
\end{description}
\end{lemma}

The symmetric part of a Banach space behaves well under contractive projections (see \cite[5.2,5.3]{Arazysurvey}).

\begin{theorem}[Stacho \cite{Stacho82}]\label{thm:stacho}
If $P$ is a contractive projection on a Banach space $X$ and $h$ is a complete holomorphic vector field on $X_0$, then $P\circ h|_{P(X)_0}$ is a complete holomorphic vector field on $P(X)_0$.  In addition $P(S(X))\subset S(P(X))$ and the partial triple product on $P(S(X))$ is given by $\{x,y,z\}=P\{x,y,z\}$ for $x,z\in P(X)$ and $y\in P(S(X))$.
\end{theorem}

Some examples of the symmetric part $S(X)$ of a Banach space $X$ are given in the seminal paper \cite{BraKauUpm78}.

\begin{itemize}
\item $X=L_p(\Omega,\Sigma,\mu)$, $1\le p<\infty$, $p\ne 2$;  $S(X)=0$
\item $X=\hbox{ (classical) }H_p$, $1\le p<\infty$, $p\ne 2$;  $S(X)=0$
\item $X=H_\infty$ (classical) or the disk algebra;  $S(X)={\bf C}$
\item $X=$ a uniform algebra $A\subset C(K)$; $S(A)=A\cap \overline{A}$
\end{itemize}

The first example above suggests Problem~\ref{prob:2} at the end of this paper.
The last example is a commutative predecessor of the example of Arazy and Solel quoted above
(\cite[Cor. 2.9(i)]{AraSol90}).
More examples, due primarily to Stacho \cite{Stacho82}, and involving Reinhardt domains are recited in \cite{Arazysurvey}, along with the following (previously) unpublished example due to Vigu\'e, showing that the symmetric part need not be complemented.

\begin{proposition}
There exists an equivalent norm on $\ell_\infty$ so that  $\ell_\infty$ in this norm has symmetric part equal to $c_0$
\end{proposition}

\section{Completely symmetric part of an operator space}

Let $A\subset B(H)$ be an operator space.  We let $K$ denote the compact operators on a separable infinite dimensional Hilbert space, say $\ell_2$. Then $K=\overline{\cup_{n=1}^\infty M_n(\CC)}$ and thus
\[
\kbta=\overline{\cup_{n=1}^\infty M_n\otimes A}=
\overline{\cup_{n=1}^\infty M_n(A)}
\]
By an abuse of notation, we shall use $\kta$ to denote 
$\cup_{n=1}^\infty M_n(A)$. We tacitly assume the embeddings $M_n(A)\subset M_{n+1}(A)\subset \kbta$ induced by adding zeros.
\smallskip

  The {\it completely symmetric part} of $A$ is defined by $CS(A)=A\cap S(K\overline{\otimes}A)$.  More precisely, if $\psi:A\rightarrow M_1(A)$ denotes the complete isometry identification, then $CS(A)=\psi^{-1}(\psi(A)\cap S(K\overline{\otimes}A))$.
\smallskip

For $1\le m<N$ let $\psi_{1,m}^N:M_N(A)\rightarrow M_N(A)$ and $\psi_{2,m}^N:M_N(A)\rightarrow M_N(A)$ be the isometries of order two defined by 
\[
\psi_{j,m}^N:\left[\begin{array}{cc}
M_m(A)&M_{m,N-m}(A)\\
M_{N-m,m}(A)&M_{N-m}(A)   
\end{array}\right]\rightarrow                                                           
\left[\begin{array}{cc}
M_m(A)&M_{m,N-m}(A)\\
M_{N-m,m}(A)&M_{N-m}(A)   
\end{array}\right]
\]
and
\[
\psi_{1,m}^N:\left[\begin{array}{cc}
a&b\\
c&d   
\end{array}\right]\rightarrow                                                           
\left[\begin{array}{cc}
a&-b\\
-c&d   
\end{array}\right]
\]
and
\[
\psi_{2,m}^N:\left[\begin{array}{cc}
a&b\\
c&d   
\end{array}\right]\rightarrow                                                           
\left[\begin{array}{cc}
a&-b\\
c&-d   
\end{array}\right].
\]

These two isometries give rise in an obvious way to two isometries $\tilde\psi_{1,m}$ and $\tilde\psi_{2,m}$ on $\kta$, which extend  to isometries $\psi_{1,m},\psi_{2,m}$ of $\kbta$ onto itself, of order 2 and fixing elementwise $M_m(A)$.  The same analysis applies to the isometries defined by, for example,
\[
\left[\begin{array}{cc}
a&b\\
c&d   
\end{array}\right]\rightarrow                                                           
\left[\begin{array}{cc}
a&b\\
-c&-d   
\end{array}\right],\left[\begin{array}{cc}
-a&-b\\
c&d   
\end{array}\right],\left[\begin{array}{cc}
-a&b\\
c&-d   
\end{array}\right].
\]

We then can define a projection $\tilde P_m$ on $\kta$ with range $M_m(A)$ via
\[
\tilde P_mx=\frac{\tilde\psi_{2,m}\left(\frac{\tilde \psi_{1,m}(x)+x}{2}\right)+\frac{\tilde\psi_{1,m}(x)+x}{2}}{2}.
\]

The  projection $\tilde P_m$ on $\kta$ extends to a projection $P_m$ on $\kbta$, with range $M_m(A)$ given by 
\[
 P_mx=\frac{\psi_{2,m}\left(\frac{\psi_{1,m}(x)+x}{2}\right)+\frac{\psi_{1,m}(x)+x}{2}}{2},
\]
or
\[
 P_m=\frac{1}{4}(\psi_{2,m}\psi_{1,m}+\psi_{2,m}+\psi_{1,m}+\hbox{Id}).
\]

\begin{proposition}\label{prop:2.1}
With the above notation,
\begin{description}
\item[(a)] $P_n(S(\kbta))=M_n(CS(A))$
\item[(b)]
 $M_n(CS(A))$ is a JB*-subtriple of $S(\kbta)$, that is\footnote{note that in the first displayed formula of (b), the triple product is the one on the JB*-triple $M_n(CS(A))$, namely,
$\tp{x}{y}{z}_{M_n(CS(A))}=P_n(\tp{x}{y}{z}_{S(\kbta)})$, which, it turns out, is actually the restriction of the triple product of $S(\kbta)$: whereas in the second displayed formula, the triple product is the partial triple product on $\kbta$},
\[
\tpc{M_n(CS(A))}{M_n(CS(A))}{M_n(CS(A))}\subset M_n(CS(A));
\]
Moreover,
\[
\tpc{M_n(A)}{M_n(CS(A))}{M_n(A)}\subset M_n(A).
\]\item[(c)] $CS(A)$ is completely isometric to  a TRO.
\end{description}
\end{proposition}
\pf\
Since $P_n$ is a linear combination of isometries of $\kbta$, and since isometries preserve the symmetric part, $P_n(S(\kbta))\subset S(\kbta)$.

Suppose $x=(x_{ij})\in P_n(S(\kbta))$. Write $x=(R_1,\cdots,R_n)^t=(C_1,\cdots,C_n)$ where $R_i,C_j$ are the rows and columns of $x$. Let $\psi_1=\psi_1^n$ and $\psi_2=\psi_2^n$ be the isometries on $\kbta$ whose action is as follows: for  $x\in M_n(A)$,
\[
\psi_1^n(x)=(R_1,-R_2,\cdots,-R_n)^t\quad,\quad \psi_2^n(x)=(-C_1,\cdots,-C_{n-1},C_n),
\]
and for an arbitrary element $y=[y_{ij}]\in\kta$, say $y\in M_N\otimes  A$, where without loss of generality $N>n$, and for $k=1,2$, $\psi_k^n$ maps $y$ into $\left[\begin{array}{cc}
\psi_k^n[y_{ij}]_{n\times n}&0\\
0&[y_{ij}]_{n<i,j\le N-n}  
\end{array}\right]$.

Then $x_{1n}\otimes e_{1n}= \frac{\psi_{2}\left(\frac{\psi_{1}(x)+x}{2}\right)+\frac{\psi_{1}(x)+x}{2}}{2}\in S(\kbta)$.

Now consider the isometry $\psi_3$ given by $\psi_3(C_1,\cdots,C_n)=(C_n,C_2,\cdots,C_{n-1},C_1)$. Then $x_{1,n}\otimes e_{11}=\psi_3(x_{1n}\otimes e_{1n})\in S(\kbta))$, and by definition, $x_{1n}\in CS(A)$.  
Continuing in this way, one sees that each $x_{ij}\in CS(A)$, proving that $P_n(S(\kbta))\subset M_n(CS(A))$

Conversely, suppose that  $x=(x_{ij})\in M_n(CS(A))$.  Since each $x_{ij}\in CS(A)$, then by definition, $x_{ij}\otimes e_{11}\in S(\kbta)$.  By using isometries as in the first part of the proof, it follows that $x_{ij}\otimes e_{ij}\in S(\kbta)$, and $x=\sum_{i,j}x_{ij}\otimes e_{ij}\in S(\kbta)$. This proves (a).
\smallskip

As noted above, $P_n$ is a contractive projection on the JB*-triple $S(\kbta)$, so that by Theorem~\ref{thm:proj}, the range of $P_n$, namely $M_n(CS(A))$, is a JB*-triple with triple product
\[
\tp{x}{y}{z}_{M_n(CS(A))}=P_n(\tp{x}{y}{z}_{S(\kbta)}),
\]
for $x,y,z\in M_n(CS(A))$.  This proves (c) by Theorem~\ref{thm:tro}.

However, $P_n$ is a linear combination of isometries of $\kbta$ which fix $M_n(A)$ elementwise, and any isometry $\psi$ of $\kbta$  preserves the partial triple product: $\psi\tp{a}{b}{c}=\tp{\psi(a)}{\psi(b)}{\psi(c)}$ for $a,c\in \kbta$ and $b\in S(\kbta)$.  This shows that 
\[
\tp{x}{y}{z}_{M_n(CS(A))}=\tp{x}{y}{z}_{S(\kbta)}
\]
for $x,y,z\in M_n(CS(A))$, proving the first part of (b). To prove the second part of (b), just note that if $x,z\in M_n(A)$ and $y\in M_n(CS(A))$, then $P_n$ fixes $\tp{x}{y}{z}$.\qed
\smallskip

\begin{corollary}
$CS(A)=M_1(CS(A))=P_1(S(\kbta))$
\end{corollary}

\begin{corollary}
$CS(A)\subset S(A)$ and $P_n\tp{y}{x}{y}=\tp{y}{x}{y}$ for $x\in M_n(CS(A))$ and $y\in M_n(A)$.
\end{corollary}
\pf\
For $x\in CS(A)$, let $\tilde x=x\otimes e_{11}$.  Then $\tilde x\in S(\kbta)$ and so there exists a complete holomorphic vector field $h_{\tilde x}$ on $(\kbta)_0$.  Since $P_1$ is a contractive projection of $\kbta$ onto $A$, by Theorem~\ref{thm:stacho}, $P_1\circ h_{\tilde x}|_{A_0}$ is a complete holomorphic vector field on $A_0$.   But $P_1\circ h_{\tilde x}|_{A_0}(0)=P_1\circ h_{\tilde x}(0)=P_1(\tilde x)=x$, proving that $x\in S(A)$.

Recall from the proof of the second part of (b) that if $x,z\in M_n(A)$ and $y\in M_n(CS(A))$, then $P_n$ fixes $\tp{x}{y}{z}$. \qed
\medskip

The symmetric part of a JC$^*$-triple coincides with the triple.
The Cartan factors of type 1 are TROs, which we have already observed are equal to their completely symmetric parts.
Since the Cartan factors of types 2,3, and 4 are not TROs it is natural to expect that their completely symmetric parts are zero.  We can verify this for finite dimensional Cartan factors of types 2, 3 and 4, as follows.

It is known \cite{Herves85} that the surjective linear isometries of the Cartan factors of types 2 and 3 are given by multiplication on the left and right by a unitary operator, and hence they are complete isometries.  The same is true for finite dimensional Cartan factors of type 4 by \cite{Vesentini92}. Using these facts and the fact that the set of inner automorphisms (hence isometries)  of any Cartan factor acts transitively on the set of  minimal tripotents (and hence on  finite rank tripotents of the same rank, \cite{HugMac09}), one can show that if  the completely symmetric part of a finite dimensional Cartan factor of type 2, 3  or 4 is not zero, then it must contain any grid which spans the Cartan factor and hence is completely isometric to a TRO.  (See Problems~\ref{prob:3} and ~\ref{prob:4} at the end of this paper.)

\section{Definition of the algebra product}

\begin{remark}
From this point on, we shall tacitly assume  the first hypothesis in our main theorem, namely that $A$ is an operator space and $v\in CS(A)$ is an element of norm 1 which satisfies  $\tp{x}{v}{v}  =x$ for every $x\in A$. With this assumption alone, we are able to construct and develop properties of the a binary product. It is not until the last step in the proof of Theorem~\ref{thm:main2} that we need to invoke the second hypothesis.

 If not explicitly stated, $a,b,c,d,x,y,z$ denote arbitrary elements of $A$.  In what follows, we work almost exclusively with $M_2(A)$, which it turns out will be sufficient for our result.  (Warning: We occasionally use capital letters $A,B,C,D$ to denote an element of the operator space $A$, as well as matrices with entries from $A$. No confusion should result as the meaning will always be clear from the context.)
\end{remark}

\begin{lemma}\label{lem:2.5}
$\tp{\left[\begin{array}{cr}
x&\pm x\\
0&0
\end{array}
\right]}{\left[\begin{array}{cr}
v&\pm v\\
0&0
\end{array}
\right]}{\left[\begin{array}{cr}
x&\pm x\\
0&0
\end{array}
\right]}=2\left[\begin{array}{cc}
\tp{x}{v}{x}&\pm \tp{x}{v}{x}\\
0&0
\end{array}
\right]$
\end{lemma}
\pf\
Let $X=\kbta$ and consider projections $Q_1$ and $Q_2$ on $X$ defined by 
$Q_1=P_{11}P_2$, $Q_2=SRP_2$ where $P_{11}$ maps $$\left[\begin{array}{cc}
a&b\\
c&d
\end{array}
\right]\hbox{ to }\left[\begin{array}{cc}
a&0\\
0&0
\end{array}
\right],$$ $S$ maps 
$$\left[\begin{array}{cc}
a&b\\
c&d
\end{array}
\right]\hbox{ to }\left[\begin{array}{cc}
a&b\\
0&0
\end{array}
\right],$$
 and $R$ maps $$ \left[\begin{array}{cc}
a&b\\
c&d
\end{array}
\right]\hbox{ to }\frac{1}{2}\left[\begin{array}{cc}
a+b&a+b\\
c+d&c+d
\end{array}
\right].$$

Let $A'=\{\left[\begin{array}{cc}
a&0\\
0&0
\end{array}
\right]:a\in A\}=Q_1X$ and $A''
=\{\left[\begin{array}{cc}
a&a\\
0&0
\end{array}
\right]:a\in A\}=Q_2X$, and let $\psi:A'\rightarrow A''$ be the isometry defined by
$\left[\begin{array}{cc}
a&0\\
0&0
\end{array}
\right]\mapsto \left[\begin{array}{cc}
a/\sqrt 2&a/\sqrt 2\\
0&0
\end{array}
\right]$.
Finally, let $v'=\left[\begin{array}{cc}
v&0\\
0&0
\end{array}
\right]$ and $v''=\left[\begin{array}{cc}
v/\sqrt 2&v/\sqrt 2\\
0&0
\end{array}
\right]$, and more generally $a'=\left[\begin{array}{cc}
a&0\\
0&0
\end{array}
\right]$, $ a''=\psi(a')=\left[\begin{array}{cc}
a/\sqrt 2&a/\sqrt 2\\
0&0
\end{array}
\right]$.

Since a surjective isometry preserves partial triple products (Lemma~\ref{lem:isometry}) and  the partial triple product on the range of a contractive projection is equal to the projection acting on the partial triple product of the original space (Theorem~\ref{thm:stacho}), we have
\[
\psi\{a'v'b'\}_{Q_1X}=\{a''v''b''\}_{Q_2X}.
\]
We unravel both sides of this equation.  In the first place
\begin{eqnarray*}
 \{a'v'b'\}_{Q_1X}&=&Q_1\{a'v'b'\}_X\\
&=&P_{11}P_2\left\{\left[\begin{array}{cc}
a&0\\
0&0
\end{array}
\right],\left[\begin{array}{cc}
v&0\\
0&0
\end{array}
\right],\left[\begin{array}{cc}
b&0\\
0&0
\end{array}
\right]\right\}_X\\
&=&P_{11}P_2\left[\begin{array}{cc}
\{avb\}&0\\
0&0
\end{array}
\right]
=
\left[\begin{array}{cc}
\{avb\}&0\\
0&0
\end{array}
\right].
\end{eqnarray*}
Thus 
\[
\psi\{a'v'b'\}_{Q_1X}=\left[\begin{array}{cc}
\{avb\}/\sqrt 2&\{avb\}/\sqrt 2\\
0&0
\end{array}
\right].
\]

Next, $R$ and $S$ are convex combinations of isometries that fix the elements of the product, so that $\{a''v''b''\}_X$ is fixed by $R$ and by $S$.  Hence,
$\{a''v''b''\}_{Q_2X}=Q_2\{a''v''b''\}_X=SRP_2\{a''v''b''\}_X=\{a''v''b''\}_X$, so that
\[
\{a''v''b''\}_{Q_2X}=\tp{\left[\begin{array}{cc}
a/\sqrt 2&a/\sqrt 2\\
0&0
\end{array}
\right]}{\left[\begin{array}{cc}
v/\sqrt 2&v/\sqrt 2\\
0&0
\end{array}
\right]}{
\left[\begin{array}{cc}
b/\sqrt 2&b/\sqrt 2\\
0&0
\end{array}
\right]}.
\]

This proves the lemma in the case of the plus sign.
The proof in the remaining case is identical, with $R$ replaced by 
$$ \left[\begin{array}{cc}
a&b\\
c&d
\end{array}
\right]\mapsto\frac{1}{2}\left[\begin{array}{cc}
a-b&b-a\\
c-d&d-c
\end{array}
\right],$$
$A''$ replaced by $
\{\left[\begin{array}{cc}
a&-a\\
0&0
\end{array}
\right]:a\in A\}$, and  $\psi$ replaced by
$\left[\begin{array}{cc}
a&0\\
0&0
\end{array}
\right]\mapsto \left[\begin{array}{cc}
a/\sqrt 2&-a/\sqrt 2\\
0&0
\end{array}
\right]$.
\qed

\begin{lemma}\label{lem:2.6}
\begin{eqnarray*}
\left[\begin{array}{cc}
\tp{x}{v}{x}&0\\
0&0
\end{array}
\right]&=&
\tp{\left[\begin{array}{cc}
0& x\\
0&0
\end{array}
\right]}{\left[\begin{array}{cc}
v&0\\
0&0
\end{array}
\right]}{\left[\begin{array}{cc}
0& x\\
0&0
\end{array}
\right]}
\\
&+&
2\tp{\left[\begin{array}{cc}
x&0\\
0&0
\end{array}
\right]}{\left[\begin{array}{cc}
0&v\\
0&0
\end{array}
\right]}{\left[\begin{array}{cc}
0&x\\
0&0
\end{array}
\right]}
\end{eqnarray*}
\end{lemma}
\pf\
By Lemma~\ref{lem:2.5}
\begin{eqnarray*}
4
\left[\begin{array}{cc}
\tp{x}{v}{x}&0\\
0&0
\end{array}
\right]
&=&
2\left[\begin{array}{cc}
\tp{x}{v}{x}& \tp{x}{v}{x}\\
0&0
\end{array}
\right]
+
2\left[\begin{array}{cc}
\tp{x}{v}{x}&- \tp{x}{v}{x}\\
0&0
\end{array}
\right]\\
&=&
\tp{\left[\begin{array}{cr}
x& x\\
0&0
\end{array}
\right]}{\left[\begin{array}{cr}
v&v\\
0&0
\end{array}
\right]}{\left[\begin{array}{cr}
x& x\\
0&0
\end{array}
\right]}\\
&+&
\tp{\left[\begin{array}{cr}
x&- x\\
0&0
\end{array}
\right]}{\left[\begin{array}{cr}
v&-v\\
0&0
\end{array}
\right]}{\left[\begin{array}{cr}
x&- x\\
0&0
\end{array}
\right]}.
\end{eqnarray*}

By expanding the right hand side of the last equation, one obtains 16 terms of which 8 cancel in pairs.  
The surviving 8 terms are

$$
2\tp{\left[\begin{array}{cr}
x& 0\\
0&0
\end{array}
\right]}{\left[\begin{array}{cr}
v&0\\
0&0
\end{array}
\right]}{\left[\begin{array}{cr}
x& 0\\
0&0
\end{array}
\right]},
$$
$$
2\tp{\left[\begin{array}{cr}
x&0\\
0&0
\end{array}
\right]}{\left[\begin{array}{cr}
0&v\\
0&0
\end{array}
\right]}{\left[\begin{array}{cr}
0&x\\
0&0
\end{array}
\right]}
$$
and
$$
4\tp{\left[\begin{array}{cr}
0&x\\
0&0
\end{array}
\right]}{\left[\begin{array}{cr}
v&0\\
0&0
\end{array}
\right]}{\left[\begin{array}{cr}
0&x\\
0&0
\end{array}
\right]}.
$$
\smallskip

Since the first surviving term above  is equal to  $2\left[\begin{array}{cr}
\tp{x}{v}{x}& 0\\
0&0
\end{array}
\right]$, the lemma is proved. \qed
\medskip

The following two lemmas, and their proofs parallel the previous two lemmas.

\begin{lemma}\label{lem:2.5prime}
$\tp{\left[\begin{array}{cr}
a&0\\
0&\pm a
\end{array}
\right]}{\left[\begin{array}{cr}
v&0\\
0&\pm v
\end{array}
\right]}{\left[\begin{array}{cr}
b&0\\
0&\pm b
\end{array}
\right]}=\left[\begin{array}{cc}
\tp{a}{v}{b}&0\\
0&\pm \tp{a}{v}{b}
\end{array}
\right]$
\end{lemma}
\pf\
Let $X=\kbta$ and consider projections $Q_1$ and $Q_2$ on $X$ defined by 
$Q_1=P_{11}P_2$, $Q_2=SRP_2$ where $P_{11}$ maps $$\left[\begin{array}{cc}
a&b\\
c&d
\end{array}
\right]\hbox{ to }\left[\begin{array}{cc}
a&0\\
0&0
\end{array}
\right],$$ $S$ maps 
$$\left[\begin{array}{cc}
a&b\\
c&d
\end{array}
\right]\hbox{ to }\left[\begin{array}{cc}
a&0\\
0&d
\end{array}
\right],$$
 and $R$ maps $$ \left[\begin{array}{cc}
a&b\\
c&d
\end{array}
\right]\hbox{ to }\frac{1}{2}\left[\begin{array}{cc}
a+d&b+c\\
b+c&a+d
\end{array}
\right].$$

Let $A'=\{\left[\begin{array}{cc}
a&0\\
0&0
\end{array}
\right]:a\in A\}=Q_1X$ and $A''
=\{\left[\begin{array}{cc}
a&0\\
0&a
\end{array}
\right]:a\in A\}=Q_2X$, and let $\psi:A'\rightarrow A''$ be the isometry defined by
$\left[\begin{array}{cc}
a&0\\
0&0
\end{array}
\right]\mapsto \left[\begin{array}{cc}
a&0\\
0&a
\end{array}
\right]$.
Finally, let $v'=\left[\begin{array}{cc}
v&0\\
0&0
\end{array}
\right]$ and $v''=\left[\begin{array}{cc}
v&0\\
0&v
\end{array}
\right]$, and more generally $a'=\left[\begin{array}{cc}
a&0\\
0&0
\end{array}
\right]$, $ a''=\psi(a')=\left[\begin{array}{cc}
a&0\\
0&a
\end{array}
\right]$.

Again by Lemma~\ref{lem:isometry} and Theorem~\ref{thm:stacho},  we have
\[
\psi\{a'v'b'\}_{Q_1X}=\{a''v''b''\}_{Q_2X}.
\]
We unravel both sides of this equation.  In the first place
\begin{eqnarray*}
 \{a'v'b'\}_{Q_1X}&=&Q_1\{a'v'b'\}_X\\
&=&P_{11}P_2\tp{\left[\begin{array}{cc}
a&0\\
0&0
\end{array}
\right]}{\left[\begin{array}{cc}
v&0\\
0&0
\end{array}
\right]}{\left[\begin{array}{cc}
b&0\\
0&0
\end{array}
\right]}_X\\
&=&P_{11}P_2\left[\begin{array}{cc}
\{avb\}&0\\
0&0
\end{array}
\right]=
\left[\begin{array}{cc}
\{avb\}&0\\
0&0
\end{array}
\right].
\end{eqnarray*}
Thus 
\[
\psi\{a'v'b'\}_{Q_1X}=\left[\begin{array}{cc}
\{avb\}&0\\
0&\{avb\}
\end{array}
\right].
\]

Next, by using appropriate isometries, for example,
$$ \left[\begin{array}{cc}
a&b\\
c&d
\end{array}
\right]\hbox{ to }\left[\begin{array}{cc}
d&b\\
c&a
\end{array}
\right].$$
  $\{a''v''b''\}_X$ is fixed by $R$ and by $S$.  Hence,
$\{a''v''b''\}_{Q_2X}=Q_2\{a''v''b''\}_X=SRP_2\{a''v''b''\}_X=\{a''v''b''\}_X$, so that
\[
\{a''v''b''\}_{Q_2X}=\tp{\left[\begin{array}{cc}
a&0\\
0&a
\end{array}
\right]}{\left[\begin{array}{cc}
v&0\\
0&v
\end{array}
\right]}{
\left[\begin{array}{cc}
b&0\\
0&b
\end{array}
\right]}.
\]

This proves the lemma in the case of the plus sign.
The proof in the remaining case is identical, with $R$ replaced by 
$$ \left[\begin{array}{cc}
a&b\\
c&d
\end{array}
\right]\mapsto\frac{1}{2}\left[\begin{array}{cc}
a-d&b-c\\
b-c&a-d
\end{array}
\right],$$
$A''$ replaced by $
\{\left[\begin{array}{cc}
a&0\\
0&-a
\end{array}
\right]:a\in A\}$, and  $\psi$ replaced by
$\left[\begin{array}{cc}
a&0\\
0&0
\end{array}
\right]\mapsto \left[\begin{array}{cc}
a&0\\
0&-a
\end{array}
\right]$.
\qed

\begin{lemma}\label{lem:2.6prime}
\[
\tp{\left[\begin{array}{cc}
x& 0\\
0&0
\end{array}
\right]}{\left[\begin{array}{cc}
0&0\\
0&v
\end{array}
\right]}{\left[\begin{array}{cc}
0& 0\\
0&y
\end{array}
\right]}
+
\tp{\left[\begin{array}{cc}
y& 0\\
0&0
\end{array}
\right]}{\left[\begin{array}{cc}
0&0\\
0&v
\end{array}
\right]}{\left[\begin{array}{cc}
0& 0\\
0&x
\end{array}
\right]}=0
\]
and
\[
\tp{\left[\begin{array}{cc}
0& 0\\
0&x
\end{array}
\right]}{\left[\begin{array}{cc}
v&0\\
0&0
\end{array}
\right]}{\left[\begin{array}{cc}
0& 0\\
0&y
\end{array}
\right]}=0.
\]
\end{lemma}
\pf\
By Lemma~\ref{lem:2.5prime}
\begin{eqnarray*}
2
\left[\begin{array}{cc}
\tp{x}{v}{x}&0\\
0&0
\end{array}
\right]
&=&
\left[\begin{array}{cc}
\tp{x}{v}{x}& 0\\
0&\tp{x}{v}{x}
\end{array}
\right]
+
\left[\begin{array}{cc}
\tp{x}{v}{x}&0\\
0&- \tp{x}{v}{x}
\end{array}
\right]\\
&=&
\tp{\left[\begin{array}{cr}
x& 0\\
0&x
\end{array}
\right]}{\left[\begin{array}{cr}
v&0\\
0&v
\end{array}
\right]}{\left[\begin{array}{cr}
x& 0\\
0&x
\end{array}
\right]}\\
&+&
\tp{\left[\begin{array}{cr}
x&0\\
0&-x
\end{array}
\right]}{\left[\begin{array}{cr}
v&0\\
0&-v
\end{array}
\right]}{\left[\begin{array}{cr}
x&0\\
0&-x
\end{array}
\right]}.
\end{eqnarray*}

By expanding the right hand side of the last equation, one obtains 16 terms of which 8 cancel in pairs.  
The surviving 8 terms are

$$
2\tp{\left[\begin{array}{cr}
x& 0\\
0&0
\end{array}
\right]}{\left[\begin{array}{cr}
v&0\\
0&0
\end{array}
\right]}{\left[\begin{array}{cr}
x& 0\\
0&0
\end{array}
\right]},
$$
$$
4\tp{\left[\begin{array}{cr}
x&0\\
0&0
\end{array}
\right]}{\left[\begin{array}{cr}
0&0\\
0&v
\end{array}
\right]}{\left[\begin{array}{cr}
0&0\\
0&x
\end{array}
\right]}
$$
and
$$
2\tp{\left[\begin{array}{cr}
0&0\\
0&x
\end{array}
\right]}{\left[\begin{array}{cr}
v&0\\
0&0
\end{array}
\right]}{\left[\begin{array}{cr}
0&0\\
0&x
\end{array}
\right]}.
$$
\smallskip

Since the first surviving term above  is equal to  $2\left[\begin{array}{cr}
\tp{x}{v}{x}& 0\\
0&0
\end{array}
\right]$, we have
\[
2\tp{\left[\begin{array}{cr}
x&0\\
0&0
\end{array}
\right]}{\left[\begin{array}{cr}
0&0\\
0&v
\end{array}
\right]}{\left[\begin{array}{cr}
0&0\\
0&x
\end{array}
\right]}
=-\tp{\left[\begin{array}{cr}
0&0\\
0&x
\end{array}
\right]}{\left[\begin{array}{cr}
v&0\\
0&0
\end{array}
\right]}{\left[\begin{array}{cr}
0&0\\
0&x
\end{array}
\right]}.
\]

Replacing $x$ by $x+y$ in this last equation results in
\[
\tp{\left[\begin{array}{cr}
x&0\\
0&0
\end{array}
\right]}{\left[\begin{array}{cr}
0&0\\
0&v
\end{array}
\right]}{\left[\begin{array}{cr}
0&0\\
0&y
\end{array}
\right]}
+\tp{\left[\begin{array}{cr}
y&0\\
0&0
\end{array}
\right]}{\left[\begin{array}{cr}
0&0\\
0&v
\end{array}
\right]}{\left[\begin{array}{cr}
0&0\\
0&x
\end{array}
\right]}
\]
\[
=
-\tp{\left[\begin{array}{cr}
0&0\\
0&x
\end{array}
\right]}{\left[\begin{array}{cr}
v&0\\
0&0
\end{array}
\right]}{\left[\begin{array}{cr}
0&0\\
0&y
\end{array}
\right]}
\]

Using the isometry of multiplication by the imaginary unit on the second row of this equation and adding then shows that both sides are zero.\qed

\begin{lemma}\label{lem:2.7}
\[
\tp{\left[\begin{array}{cc}
x& 0\\
0&0
\end{array}
\right]}{\left[\begin{array}{cc}
0&0\\
0&v
\end{array}
\right]}{\left[\begin{array}{cc}
a& b\\
c&0
\end{array}
\right]}=0
\]
and\footnote{This is true with the second $v$ replaced by an arbitrary element of $A$. We shall  prove it  in Corollary~\ref{cor:4.4} (using the special case given by Lemma~\ref{lem:2.7}) }
\[
\tp{\left[\begin{array}{cc}
x& 0\\
0&0
\end{array}
\right]}{\left[\begin{array}{cc}
0&0\\
0&v
\end{array}
\right]}{\left[\begin{array}{cc}
0& 0\\
0&v
\end{array}
\right]}=0 ,
\]
Equivalently,
\[
\tp{\left[\begin{array}{cc}
0& x\\
0&0
\end{array}
\right]}{\left[\begin{array}{cc}
0&0\\
v&0
\end{array}
\right]}{\left[\begin{array}{cc}
a& b\\
0&d
\end{array}
\right]}=0
\]
and
\[
\tp{\left[\begin{array}{cc}
0& x\\
0&0
\end{array}
\right]}{\left[\begin{array}{cc}
0&0\\
v&0
\end{array}
\right]}{\left[\begin{array}{cc}
0& 0\\
v&0
\end{array}
\right]}=0,
\]
\end{lemma}
\pf\
The second statement follows from the first  by using the isometry
\[
\left[\begin{array}{cc}
a&b\\
c&d   
\end{array}\right]\mapsto                                                           
\left[\begin{array}{cc}
b&a\\
d&c   
\end{array}\right].
\]
Using Lemma~\ref{lem:2.6prime} and an appropriate isometry (interchange both rows and columns simultaneously)  yields
\[
\tp{\left[\begin{array}{cc}
x& 0\\
0&0
\end{array}
\right]}{\left[\begin{array}{cc}
0&0\\
0&v
\end{array}
\right]}{\left[\begin{array}{cc}
a& 0\\
0&0
\end{array}
\right]}=0.
\]
Next, the isometry
\[
\left[\begin{array}{cc}
a&b\\
c&d   
\end{array}\right]\mapsto                                                           
\left[\begin{array}{cc}
-a&-b\\
c&d   
\end{array}\right].
\]
shows that
\[
\tp{\left[\begin{array}{cc}
x& 0\\
0&0
\end{array}
\right]}{\left[\begin{array}{cc}
0&0\\
0&v
\end{array}
\right]}{\left[\begin{array}{cc}
0&b\\
0&0
\end{array}
\right]}
=\matr{0}{0}{C}{D},
\]
for some $C,D \in A$.    Similarly, the isometry
\[
\left[\begin{array}{cc}
a&b\\
c&d   
\end{array}\right]\mapsto                                                           
\left[\begin{array}{cc}
a&-b\\
c&-d   
\end{array}\right]
\]
shows that 
\[
\tp{\left[\begin{array}{cc}
x& 0\\
0&0
\end{array}
\right]}{\left[\begin{array}{cc}
0&0\\
0&v
\end{array}
\right]}{\left[\begin{array}{cc}
0&b\\
0&0
\end{array}
\right]}
=\matr{0}{0}{C}{0}.
\]

Applying the isometry of multiplication of the second row by the imaginary unit shows that $C=0$.  Hence
\[
\tp{\left[\begin{array}{cc}
x& 0\\
0&0
\end{array}
\right]}{\left[\begin{array}{cc}
0&0\\
0&v
\end{array}
\right]}{\left[\begin{array}{cc}
0& b\\
0&0
\end{array}
\right]}=0. 
\]
By appropriate  use of  isometries as above,
\[
\tp{\left[\begin{array}{cc}
x& 0\\
0&0
\end{array}
\right]}{\left[\begin{array}{cc}
0&0\\
0&v
\end{array}
\right]}{\left[\begin{array}{cc}
0&0\\
c&0
\end{array}
\right]}
=\matr{0}{B}{0}{0}
\]
for some $B\in A$.
Applying the isometry of multiplication of the second column by the imaginary unit  shows that $B=0$.  Hence
\[\tp{\left[\begin{array}{cc}
x& 0\\
0&0
\end{array}
\right]}{\left[\begin{array}{cc}
0&0\\
0&v
\end{array}
\right]}{\left[\begin{array}{cc}
0& 0\\
c&0
\end{array}
\right]}=0.
\]
It remains to show that
\[
\tp{\left[\begin{array}{cc}
x& 0\\
0&0
\end{array}
\right]}{\left[\begin{array}{cc}
0&0\\
0&v
\end{array}
\right]}{\left[\begin{array}{cc}
0& 0\\
0&v
\end{array}
\right]}=0,
\]
To this end, by the main identity,
\begin{equation}\label{eq:1214111}
\tp 
{\matr{v}{0}{0}{0}}   
{\matr{v}{0}{0}{0}}   
{\tp
{\matr{v}{0}{0}{0}}
{\matr{0}{0}{0}{v}}
{\matr{0}{0}{0}{x}}
}
=R-S+T
\end{equation}
where
\[
R=\tp
{\tp
{\matr{v}{0}{0}{0}}
{\matr{v}{0}{0}{0}}
{\matr{v}{0}{0}{0}}
}
{\matr{0}{0}{0}{v}}
{\matr{0}{0}{0}{x}
},
\]
\[
S=\tp
{\matr{v}{0}{0}{0}}
{\tp
{\matr{v}{0}{0}{0}}
{\matr{v}{0}{0}{0}}
{\matr{0}{0}{0}{v}}
}
{\matr{0}{0}{0}{x}}
\]
and
\[
T=\tp
{\matr{v}{0}{0}{0}}
{\matr{0}{0}{0}{v}}
{
\tp
{\matr{v}{0}{0}{0}}
{\matr{v}{0}{0}{0}}
{\matr{0}{0}{0}{x}}
}.
\]

Since 
\begin{equation}\label{eq:1220111}
\tp
{\matr{v}{0}{0}{0}}
{\matr{0}{0}{0}{v}}
{\matr{0}{0}{0}{x}}
\hbox{
has the form } \matr{A}{0}{0}{0},
\end{equation}
 the left side of (\ref{eq:1214111}) is equal to 
\begin{equation}\label{eq:1220112}
\tp
{\matr{v}{0}{0}{0}}
{\matr{v}{0}{0}{0}}
{\matr{A}{0}{0}{0}}
=\matr{\tp{v}{v}{A}}{0}{0}{0}=\matr{A}{0}{0}{0}.
\end{equation}
This term is also equal to $R$ since 
\[
\tp
{\matr{v}{0}{0}{0}}
{\matr{v}{0}{0}{0}}
{\matr{v}{0}{0}{0}}=\matr{v}{0}{0}{0}.
\]
Since $S=0$, we have $T=0$.
We next apply the main identity  to get 
$0=T=R'-S'+T'$, where

\[
R'=\tp
{\tp
{\matr{v}{0}{0}{0}}
{\matr{0}{v}{0}{0}}
{\matr{v}{0}{0}{0}}
}
{\matr{v}{0}{0}{0}}
{\matr{0}{0}{0}{x}
},
\]
\[
S'=\tp
{\matr{v}{0}{0}{0}}
{\tp
{\matr{0}{0}{0}{v}}
{\matr{v}{0}{0}{0}}
{\matr{v}{0}{0}{0}}
}
{\matr{0}{0}{0}{x}}
\]
and
\[
T'=\tp
{\matr{v}{0}{0}{0}}
{\matr{v}{0}{0}{0}}
{
\tp
{\matr{v}{0}{0}{0}}
{\matr{0}{0}{0}{v}}
{\matr{0}{0}{0}{x}}
}.
\]
By direct calculation, $R'=0$ and $S'=0$, and since $T=0$ we have  $T'=0$ so that  by (\ref{eq:1220111}) and (\ref{eq:1220112}),
\begin{eqnarray*}
0&=&\tp
{\matr{v}{0}{0}{0}}
{\matr{v}{0}{0}{0}}
{
\tp
{\matr{v}{0}{0}{0}}
{\matr{0}{0}{0}{v}}
{\matr{0}{0}{0}{x}}
}\\
&=&\tp
{\matr{v}{0}{0}{0}}
{\matr{0}{0}{0}{v}}
{\matr{0}{0}{0}{x}}.\qed
\end{eqnarray*}

\begin{lemma}\label{lem:2.8}
\begin{eqnarray*}
\left[\begin{array}{cc}
\tp{x}{v}{y}&0\\
0&0
\end{array}
\right]&=&
\tp{\left[\begin{array}{cc}
0& x\\
0&0
\end{array}
\right]}{\left[\begin{array}{cc}
v&0\\
0&0
\end{array}
\right]}{\left[\begin{array}{cc}
0& y\\
0&0
\end{array}
\right]}
\\
&+&
\tp{\left[\begin{array}{cc}
x&0\\
0&0
\end{array}
\right]}{\left[\begin{array}{cc}
0&v\\
0&0
\end{array}
\right]}{\left[\begin{array}{cc}
0&y\\
0&0
\end{array}
\right]}\\
&+&
\tp{\left[\begin{array}{cc}
0&x\\
0&0
\end{array}
\right]}{\left[\begin{array}{cc}
0&v\\
0&0
\end{array}
\right]}{\left[\begin{array}{cc}
y&0\\
0&0
\end{array}
\right]}
\end{eqnarray*}
\end{lemma}
\pf\
Replace $x$ in Lemma~\ref{lem:2.6} by $x+y$.\qed
\smallskip

We can repeat some of the preceding arguments to obtain the following three lemmas, which will be used in the proof of 
Lemma~\ref{lem:2.13}.
The proof of the following lemma is, except for notation,
identical to those of Lemma~\ref{lem:2.5} and Lemma~\ref{lem:2.5prime}.

\begin{lemma}\label{lem:2.10}
$\tp{\left[\begin{array}{cc}
0& x\\
0&\pm x
\end{array}
\right]}{\left[\begin{array}{cc}
0&v\\
0&\pm v
\end{array}
\right]}{\left[\begin{array}{cc}
0&x\\
0&\pm x
\end{array}
\right]}=\left[\begin{array}{cc}
0&2 \tp{x}{v}{x}\\
0&\pm 2\tp{x}{v}{x}
\end{array}
\right]$
\end{lemma}

The proof of the following lemma parallels exactly the proof of Lemma~\ref{lem:2.6}, using Lemma~\ref{lem:2.10} in place of Lemma~\ref{lem:2.5}.

\begin{lemma}\label{lem:2.11}
\begin{eqnarray*}
\left[\begin{array}{cc}
0&\tp{x}{v}{x}\\
0&0
\end{array}
\right]&=&
\tp{\left[\begin{array}{cc}
0&0 \\
0&x
\end{array}
\right]}{\left[\begin{array}{cc}
0&v\\
0&0
\end{array}
\right]}{\left[\begin{array}{cc}
0&0 \\
0&x
\end{array}
\right]}
\\
&+&
2\tp{\left[\begin{array}{cc}
0&0\\
0&x
\end{array}
\right]}{\left[\begin{array}{cc}
0&0\\
0&v
\end{array}
\right]}{\left[\begin{array}{cc}
0&x\\
0&0
\end{array}
\right]}
\end{eqnarray*}
\end{lemma}

As in Lemma~\ref{lem:2.8}, polarization  of Lemma~\ref{lem:2.11} yields the following lemma.

\begin{lemma}\label{lem:2.12}
\begin{eqnarray*}
\left[\begin{array}{cc}
0&\tp{x}{v}{y}\\
0&0
\end{array}
\right]&=&
\tp{\left[\begin{array}{cc}
0&0 \\
0&x
\end{array}
\right]}{\left[\begin{array}{cc}
0&v\\
0&0
\end{array}
\right]}{\left[\begin{array}{cc}
0&0 \\
0&y
\end{array}
\right]}
\\
&+&
\tp{\left[\begin{array}{cc}
0&0\\
0&x
\end{array}
\right]}{\left[\begin{array}{cc}
0&0\\
0&v
\end{array}
\right]}{\left[\begin{array}{cc}
0&y\\
0&0
\end{array}
\right]}
\\
&+&
\tp{\left[\begin{array}{cc}
0&0\\
0&y
\end{array}
\right]}{\left[\begin{array}{cc}
0&0\\
0&v
\end{array}
\right]}{\left[\begin{array}{cc}
0&x\\
0&0
\end{array}
\right]}
\end{eqnarray*}
\end{lemma}

\begin{lemma}\label{lem:2.13}
\[
\tp{\left[\begin{array}{cc}
0&v\\
0&0
\end{array}
\right]}{\left[\begin{array}{cc}
v&0\\
0&0
\end{array}
\right]}{\left[\begin{array}{cc}
0&x\\
0&0
\end{array}
\right]}=0.
\]
\end{lemma}
\pf\
Set $y=v$ in Lemma~\ref{lem:2.12} and apply $D\left(\left[\begin{array}{cc}
0&v\\
0&0
\end{array}
\right],\left[\begin{array}{cc}
v&0\\
0&0
\end{array}
\right]\right)$ to each side of the equation in that lemma.

The three terms on the right each vanish, as is seen by applying the main identity to each term and making use of  Lemma~\ref{lem:2.7}, and  the fact that $CS(A)$ is a TRO, and hence $M_2(CS(A))$ is a JB$^*$-triple.  The lemma is proved.
\smallskip

For the sake of clarity, we again include the details of the proof.
Explicitly,
\[
\tp{\left[\begin{array}{cc}
0&v\\
0&0
\end{array}
\right]}{\left[\begin{array}{cc}
v&0\\
0&0
\end{array}
\right]}{\left[\begin{array}{cc}
0&x\\
0&0
\end{array}
\right]}=R+S+T,
\]
 where
\[
R=
\tp{{\left[\begin{array}{cc}
0&v\\
0&0
\end{array}
\right]}}{{\left[\begin{array}{cc}
v&0\\
0&0
\end{array}
\right]}}
{\tp{\left[\begin{array}{cc}
0& 0\\
0&x
\end{array}
\right]}{\left[\begin{array}{cc}
0&v\\
0&0
\end{array}
\right]}{\left[\begin{array}{cc}
0&0\\
0&v
\end{array}
\right]}
}
\]
\[
S=
\tp{{\left[\begin{array}{cc}
0& v\\
0&0
\end{array}
\right]}}{{\left[\begin{array}{cc}
v&0\\
0&0
\end{array}
\right]}}
{\tp{\left[\begin{array}{cc}
0& 0\\
0&x
\end{array}
\right]}{\left[\begin{array}{cc}
0&0\\
0&v
\end{array}
\right]}{\left[\begin{array}{cc}
0&v\\
0&0
\end{array}
\right]}
}
\]
and 
\[
T=
\tp{{\left[\begin{array}{cc}
0& v\\
0&0
\end{array}
\right]}}{{\left[\begin{array}{cc}
v&0\\
0&0
\end{array}
\right]}}
{\tp{\left[\begin{array}{cc}
0& 0\\
0&v
\end{array}
\right]}{\left[\begin{array}{cc}
0&0\\
0&v
\end{array}
\right]}{\left[\begin{array}{cc}
0&x\\
0&0
\end{array}
\right]}
}.
\]

By the main identity,
\[
R=\tp{\tp{\left[\begin{array}{cc}
0&v\\
0&0
\end{array}
\right]}{\left[\begin{array}{cc}
v&0\\
0&0
\end{array}
\right]}{\left[\begin{array}{cc}
0&0\\
0&x
\end{array}
\right]}}{\left[\begin{array}{cc}
0&v\\
0&0
\end{array}
\right]}{\left[\begin{array}{cc}
0&0\\
0&v
\end{array}
\right]}
\]
\[
-\tp{\left[\begin{array}{cc}
0&0\\
0&x
\end{array}
\right]}{\tp{\left[\begin{array}{cc}
v&0\\
0&0
\end{array}
\right]}{\left[\begin{array}{cc}
0&v\\
0&0
\end{array}
\right]}{\left[\begin{array}{cc}
0&v\\
0&0
\end{array}
\right]}}{\left[\begin{array}{cc}
0&0\\
0&v
\end{array}
\right]}
\]
\[
+\tp{\left[\begin{array}{cc}
0&0\\
0&x
\end{array}
\right]}{\left[\begin{array}{cc}
0&v\\
0&0
\end{array}
\right]}{\tp{\left[\begin{array}{cc}
0&v\\
0&0
\end{array}
\right]}{\left[\begin{array}{cc}
v&0\\
0&0
\end{array}
\right]}{\left[\begin{array}{cc}
0&0\\
0&v
\end{array}
\right]}}.
\]
The first term is zero by Lemma~\ref{lem:2.7}. The third term is zero by direct calculation in $CS(A)$ as in Proposition~\ref{prop:2.1}.
 The middle term is zero by Lemma~\ref{lem:2.7} since by Proposition~\ref{prop:2.1},
\[
\tp{\left[\begin{array}{cc}
v&0\\
0&0
\end{array}
\right]}{\left[\begin{array}{cc}
0&v\\
0&0
\end{array}
\right]}{\left[\begin{array}{cc}
0&v\\
0&0
\end{array}
\right]}=\frac{1}{2}\left[\begin{array}{cc}
v&0\\
0&0
\end{array}
\right].
\]
Hence $R=0$.

Again by the main identity,
\[
S=\tp{\tp{\left[\begin{array}{cc}
0&v\\
0&0
\end{array}
\right]}{\left[\begin{array}{cc}
v&0\\
0&0
\end{array}
\right]}{\left[\begin{array}{cc}
0&0\\
0&x
\end{array}
\right]}}{\left[\begin{array}{cc}
0&0\\
0&v
\end{array}
\right]}{\left[\begin{array}{cc}
0&v\\
0&0
\end{array}
\right]}
\]
\[
-\tp{\left[\begin{array}{cc}
0&0\\
0&x
\end{array}
\right]}{\tp{\left[\begin{array}{cc}
v&0\\
0&0
\end{array}
\right]}{\left[\begin{array}{cc}
0&v\\
0&0
\end{array}
\right]}{\left[\begin{array}{cc}
0&0\\
0&v
\end{array}
\right]}}{\left[\begin{array}{cc}
0&v\\
0&0
\end{array}
\right]}
\]
\[
+\tp{\left[\begin{array}{cc}
0&0\\
0&x
\end{array}
\right]}{\left[\begin{array}{cc}
0&0\\
0&v
\end{array}
\right]}{\tp{\left[\begin{array}{cc}
0&v\\
0&0
\end{array}
\right]}{\left[\begin{array}{cc}
v&0\\
0&0
\end{array}
\right]}{\left[\begin{array}{cc}
0&v\\
0&0
\end{array}
\right]}}.
\]
The 
first term is zero by Lemma~\ref{lem:2.7} and the
third 
term is zero by direct calculation:
\[
\tp{\left[\begin{array}{cc}
0&v\\
0&0
\end{array}
\right]}{\left[\begin{array}{cc}
v&0\\
0&0
\end{array}
\right]}{\left[\begin{array}{cc}
0&v\\
0&0
\end{array}
\right]}=\left[\begin{array}{cc}
0&v\\
0&0
\end{array}
\right]\left[\begin{array}{cc}
v^*&0\\
0&0
\end{array}
\right]\left[\begin{array}{cc}
0&v\\
0&0
\end{array}
\right]=0.
\]

 The middle term is also zero by Lemma~\ref{lem:2.7} since by Proposition~\ref{prop:2.1},
\[
\tp{\left[\begin{array}{cc}
v&0\\
0&0
\end{array}
\right]}{\left[\begin{array}{cc}
0&v\\
0&0
\end{array}
\right]}{\left[\begin{array}{cc}
0&0\\
0&v
\end{array}
\right]}=\frac{1}{2}\left[\begin{array}{cc}
0&0\\
v&0
\end{array}
\right].
\]
Hence $S=0$. 

Finally, again by the main identity,
\[
T=\tp{\tp{\left[\begin{array}{cc}
0&v\\
0&0
\end{array}
\right]}{\left[\begin{array}{cc}
v&0\\
0&0
\end{array}
\right]}{\left[\begin{array}{cc}
0&0\\
0&v
\end{array}
\right]}}{\left[\begin{array}{cc}
0&0\\
0&v
\end{array}
\right]}{\left[\begin{array}{cc}
0&x\\
0&0
\end{array}
\right]}
\]
\[
-\tp{\left[\begin{array}{cc}
0&0\\
0&v
\end{array}
\right]}{\tp{\left[\begin{array}{cc}
v&0\\
0&0
\end{array}
\right]}{\left[\begin{array}{cc}
0&v\\
0&0
\end{array}
\right]}{\left[\begin{array}{cc}
0&0\\
0&v
\end{array}
\right]}}{\left[\begin{array}{cc}
0&x\\
0&0
\end{array}
\right]}
\]
\[
+\tp{\left[\begin{array}{cc}
0&0\\
0&v
\end{array}
\right]}{\left[\begin{array}{cc}
0&0\\
0&v
\end{array}
\right]}{\tp{\left[\begin{array}{cc}
0&v\\
0&0
\end{array}
\right]}{\left[\begin{array}{cc}
v&0\\
0&0
\end{array}
\right]}{\left[\begin{array}{cc}
0&x\\
0&0
\end{array}
\right]}}.
\]
The first term is zero by  direct calculation
 in $CS(A)$ as in Proposition~\ref{prop:2.1}.
  The third term is of the form
\[
\tp{\left[\begin{array}{cc}
0&0\\
0&v
\end{array}
\right]}{\left[\begin{array}{cc}
0&0\\
0&v
\end{array}
\right]}{\left[\begin{array}{cc}
A&0\\
0&0
\end{array}
\right]}
\]
so it is zero by Lemma~\ref{lem:2.7}. The middle term is zero
by Lemma~\ref{lem:2.7} since by Proposition~\ref{prop:2.1},
\[
\tp{\left[\begin{array}{cc}
v&0\\
0&0
\end{array}
\right]}{\left[\begin{array}{cc}
0&v\\
0&0
\end{array}
\right]}{\left[\begin{array}{cc}
0&0\\
0&v
\end{array}
\right]}=\frac{1}{2}\left[\begin{array}{cc}
0&0\\
v&0
\end{array}
\right].
\]
Hence $T=0$. \qed
\smallskip

\begin{lemma}\label{lem:2.14}
\[
\tp{\left[\begin{array}{cc}
0&x\\
0&0
\end{array}
\right]}{\left[\begin{array}{cc}
v&0\\
0&0
\end{array}
\right]}{\left[\begin{array}{cc}
0&y\\
0&0
\end{array}
\right]}=0.
\]
\end{lemma}
\pf\ 
By applying the isometries of multiplication of the second column and second row by $-1$, we see that
\[
\tp{\left[\begin{array}{cc}
0&x\\
0&0
\end{array}
\right]}{\left[\begin{array}{cc}
v&0\\
0&0
\end{array}
\right]}{\left[\begin{array}{cc}
0&y\\
0&0
\end{array}
\right]}=\left[\begin{array}{cc}
a&0\\
0&0
\end{array}
\right]
\]
and that
\[
\tp{\left[\begin{array}{cc}
x&0\\
0&0
\end{array}
\right]}{\left[\begin{array}{cc}
0&v\\
0&0
\end{array}
\right]}{\left[\begin{array}{cc}
y&0\\
0&0
\end{array}
\right]}=\left[\begin{array}{cc}
0&a\\
0&0
\end{array}
\right].
\]

By Lemma~\ref{lem:2.8}

\begin{eqnarray}\label{star}
\left[\begin{array}{cc}
a&0\\
0&0
\end{array}
\right]=\left[\begin{array}{cc}
\tp{a}{v}{v}&0\\
0&0
\end{array}
\right]&=&
\tp{\left[\begin{array}{cc}
0& a\\
0&0
\end{array}
\right]}{\left[\begin{array}{cc}
v&0\\
0&0
\end{array}
\right]}{\left[\begin{array}{cc}
0& v\\
0&0
\end{array}
\right]}
\\\nonumber
&+&
\tp{\left[\begin{array}{cc}
a&0\\
0&0
\end{array}
\right]}{\left[\begin{array}{cc}
0&v\\
0&0
\end{array}
\right]}{\left[\begin{array}{cc}
0&v\\
0&0
\end{array}
\right]}
\\\nonumber
&+&
\tp{\left[\begin{array}{cc}
0& a\\
0&0
\end{array}
\right]}{\left[\begin{array}{cc}
0&v\\
0&0
\end{array}
\right]}{\left[\begin{array}{cc}
v&0\\
0&0
\end{array}
\right]}
\end{eqnarray}

The first term on the right side of (\ref{star}) is zero by Lemma~\ref{lem:2.13}.

Let us write the second term on the right side of  (\ref{star}) as
\[
\tp{\left[\begin{array}{cc}
a& 0\\
0&0
\end{array}
\right]}{\left[\begin{array}{cc}
0&v\\
0&0
\end{array}
\right]}{\left[\begin{array}{cc}
0&v\\
0&0
\end{array}
\right]}=
\tp{{\left[\begin{array}{cc}
0& v\\
0&0
\end{array}
\right]}}{{\left[\begin{array}{cc}
0&v\\
0&0
\end{array}
\right]}}
{\tp{\left[\begin{array}{cc}
0& x\\
0&0
\end{array}
\right]}{\left[\begin{array}{cc}
v&0\\
0&0
\end{array}
\right]}{\left[\begin{array}{cc}
0&y\\
0&0
\end{array}
\right]}
}
\]
and apply the main Jordan identity to the right side, which we write  symbolically as $\tp{A}{B}{\tp{C}{D}{E}}$ to obtain
\begin{equation}\label{starstar}
\tp{A}{B}{\tp{C}{D}{E}}=\tp{\tp{A}{B}{C}}{D}{E}-\tp{C}{\tp{B}{A}{D}}{E}+\tp{C}{D}{\tp{A}{B}{E}}
\end{equation}
We then calculate each term on the right side of (\ref{starstar})  to obtain
\[
\tp{\tp{A}{B}{C}}{D}{E}=\tp{\left[\begin{array}{cc}
0& \tp{v}{v}{x}\\
0&0
\end{array}
\right]}{\left[\begin{array}{cc}
v&0\\
0&0
\end{array}
\right]}{\left[\begin{array}{cc}
0&y\\
0&0
\end{array}
\right]}=\tp{\left[\begin{array}{cc}
0& x\\
0&0
\end{array}
\right]}{\left[\begin{array}{cc}
v&0\\
0&0
\end{array}
\right]}{\left[\begin{array}{cc}
0&y\\
0&0
\end{array}
\right]}
\]
\[
\tp{C}{\tp{B}{A}{D}}{E}=\frac{1}{2}\tp{\left[\begin{array}{cc}
0& x\\
0&0
\end{array}
\right]}{\left[\begin{array}{cc}
vv^*v&0\\
0&0
\end{array}
\right]}{\left[\begin{array}{cc}
0&y\\
0&0
\end{array}
\right]}=\frac{1}{2}\tp{\left[\begin{array}{cc}
0& x\\
0&0
\end{array}
\right]}{\left[\begin{array}{cc}
v&0\\
0&0
\end{array}
\right]}{\left[\begin{array}{cc}
0&y\\
0&0
\end{array}
\right]}
\]
\[
\tp{C}{D}{\tp{A}{B}{E}}=\tp{\left[\begin{array}{cc}
0& x\\
0&0
\end{array}
\right]}{\left[\begin{array}{cc}
v&0\\
0&0
\end{array}
\right]}{
\tp{\left[\begin{array}{cc}
0& v\\
0&0
\end{array}
\right]}{\left[\begin{array}{cc}
0&v\\
0&0
\end{array}
\right]}{\left[\begin{array}{cc}
0&y\\
0&0
\end{array}
\right]}}
=
\tp{\left[\begin{array}{cc}
0& x\\
0&0
\end{array}
\right]}{\left[\begin{array}{cc}
v&0\\
0&0
\end{array}
\right]}{\left[\begin{array}{cc}
0&y\\
0&0
\end{array}
\right]}.
\]
The second term on the right side of  (\ref{star}) is therefore equal to $$\frac{3}{2}\tp{\left[\begin{array}{cc}
0& x\\
0&0
\end{array}
\right]}{\left[\begin{array}{cc}
v&0\\
0&0
\end{array}
\right]}{\left[\begin{array}{cc}
0&y\\
0&0
\end{array}
\right]}$$

Let us write the third  term on the right side of (\ref{star}) as
\[
\tp{\left[\begin{array}{cc}
0& a\\
0&0
\end{array}
\right]}{\left[\begin{array}{cc}
0&v\\
0&0
\end{array}
\right]}{\left[\begin{array}{cc}
v&0\\
0&0
\end{array}
\right]}=
\tp{{\left[\begin{array}{cc}
v& 0\\
0&0
\end{array}
\right]}}{{\left[\begin{array}{cc}
0&v\\
0&0
\end{array}
\right]}}
{\tp{\left[\begin{array}{cc}
x& 0\\
0&0
\end{array}
\right]}{\left[\begin{array}{cc}
0&v\\
0&0
\end{array}
\right]}{\left[\begin{array}{cc}
y&0\\
0&0
\end{array}
\right]}
}
\]
and apply the main Jordan identity to the right side, which we again write  symbolically as $\tp{A'}{B'}{\tp{C'}{D'}{E'}}$ to obtain
\[
\tp{A'}{B'}{\tp{C'}{D'}{E'}}=\tp{\tp{A'}{B'}{C'}}{D'}{E'}-\tp{C'}{\tp{B'}{A'}{D'}}{E'}+\tp{C'}{D'}{\tp{A'}{B'}{E'}}
\]
We then calculate each term on the right side and find that each of these terms vanishes, the first and third by Lemma~\ref{lem:2.13} and the second by the fact that $CS(A)$ is a TRO.

We have thus shown that
$$\tp{\left[\begin{array}{cc}
0& x\\
0&0
\end{array}
\right]}{\left[\begin{array}{cc}
v&0\\
0&0
\end{array}
\right]}{\left[\begin{array}{cc}
0&y\\
0&0
\end{array}
\right]}=\frac{3}{2}\tp{\left[\begin{array}{cc}
0& x\\
0&0
\end{array}
\right]}{\left[\begin{array}{cc}
v&0\\
0&0
\end{array}
\right]}{\left[\begin{array}{cc}
0&y\\
0&0
\end{array}
\right]},$$
proving the lemma.\qed
\smallskip

\begin{definition}\label{def:3.13}
Let us now define a product $y\cdot x$ by
\[
\left[\begin{array}{cc}
y\cdot x& 0\\
0&0
\end{array}
\right]=2
\tp{\left[\begin{array}{cc}
x& 0\\
0&0
\end{array}
\right]}{\left[\begin{array}{cc}
0&v\\
0&0
\end{array}
\right]}{\left[\begin{array}{cc}
0&y\\
0&0
\end{array}
\right]}
\]
and denote the corresponding matrix product by $X\cdot Y$. 
That is, if $X=[x_{ij}]$ and $Y=[y_{ij}]$, then $X\cdot Y=[z_{ij}]$ where
\[
z_{ij}=\sum_k x_{ik}\cdot y_{kj}.
\]
\end{definition}

 Note that 
\begin{equation}\label{eq:0606121}
\tp{x}{v}{y}=\frac{1}{2}(y\cdot x+x\cdot y).
\end{equation}
since
by Lemmas~\ref{lem:2.8} and \ref{lem:2.14}
we can write
\smallskip
\[
\left[\begin{array}{cc}
\tp{x}{v}{y}& 0\\
0&0
\end{array}
\right]=
\tp{\left[\begin{array}{cc}
x& 0\\
0&0
\end{array}
\right]}{\left[\begin{array}{cc}
0&v\\
0&0
\end{array}
\right]}{\left[\begin{array}{cc}
0&y\\
0&0
\end{array}
\right]}
+
\tp{\left[\begin{array}{cc}
0& x\\
0&0
\end{array}
\right]}{\left[\begin{array}{cc}
0&v\\
0&0
\end{array}
\right]}{\left[\begin{array}{cc}
y&0\\
0&0
\end{array}
\right]}.
\]

\section{Main result}

The following lemma, in which the right side  is equal to  $\frac{1}{2}\left[\begin{array}{cc}
0& 0\\
0&x\cdot y
\end{array}
\right]$, is needed to prove 
that $v$ is a unit element for the product $x\cdot y$ (Lemma~\ref{lem:1217111}), and to prove
Proposition~\ref{prop:2.14} below, which is another key step in the proof.

\begin{lemma}\label{lem:2.18}
\[
\tp{\left[\begin{array}{cc}
0& 0\\
x&0
\end{array}
\right]}{\left[\begin{array}{cc}
v&0\\
0&0
\end{array}
\right]}{\left[\begin{array}{cc}
0&y\\
0&0
\end{array}
\right]}
=
\tp{\left[\begin{array}{cc}
0& 0\\
x&0
\end{array}
\right]}{\left[\begin{array}{cc}
0&0\\
v&0
\end{array}
\right]}{\left[\begin{array}{cc}
0&0\\
0&y
\end{array}
\right]}
\]
\end{lemma}
\pf\ 
Let $\psi$ be  the isometry
\[
\left[\begin{array}{cc}
x& y\\
0&0
\end{array}
\right]\mapsto\frac{1}{\sqrt 2} \left[\begin{array}{cc}
x& y\\
x&y
\end{array}
\right].
\]
As in the proofs of 
Lemmas~\ref{lem:2.5},\ref{lem:2.5prime} and \ref{lem:2.10}, $\psi$ preserves partial triple products.  Thus,
\begin{eqnarray*}\frac{1}{2\sqrt 2}\matr{0}{x\cdot y}{0}{x\cdot y}&=&\frac{1}{2}\psi\left(\matr{0}{x\cdot y}{0}{0}\right)\\
&=&
\psi\left(\tp{\matr{x}{0}{0}{0}}{\matr{v}{0}{0}{0}}{\matr{0}{y}{0}{0}}\right)\\
&=&\left(\frac{1}{\sqrt 2}\right)^3\tp{\matr{x}{0}{x}{0}}{\matr{v}{0}{v}{0}}{\matr{0}{y}{0}{y}}\\
&=&\left(\frac{1}{\sqrt 2}\right)^3\left(\tp{\matr{x}{0}{0}{0}}{\matr{v}{0}{0}{0}}{\matr{0}{y}{0}{0}}\right)\\
&+&\left(\frac{1}{\sqrt 2}\right)^3\left(\tp{\matr{0}{0}{x}{0}}{\matr{v}{0}{0}{0}}{\matr{0}{y}{0}{0}}\right)\\
&+&\left(\frac{1}{\sqrt 2}\right)^3\left(\tp{\matr{x}{0}{0}{0}}{\matr{0}{0}{v}{0}}{\matr{0}{0}{0}{y}}\right)\\
&+&\left(\frac{1}{\sqrt 2}\right)^3\left(\tp{\matr{0}{0}{x}{0}}{\matr{0}{0}{v}{0}}{\matr{0}{0}{0}{y}}\right),
\end{eqnarray*}
so that
\begin{eqnarray*}
\matr{0}{x\cdot y}{0}{x\cdot y}
&=&\tp{\matr{x}{0}{0}{0}}{\matr{v}{0}{0}{0}}{\matr{0}{y}{0}{0}}\\
&+&\tp{\matr{0}{0}{x}{0}}{\matr{v}{0}{0}{0}}{\matr{0}{y}{0}{0}}\\
&+&\tp{\matr{x}{0}{0}{0}}{\matr{0}{0}{v}{0}}{\matr{0}{0}{0}{y}}\\
&+&(\tp{\matr{0}{0}{x}{0}}{\matr{0}{0}{v}{0}}{\matr{0}{0}{0}{y}}.
\end{eqnarray*}
On the other hand,
\begin{eqnarray*}
\matr{0}{x\cdot y}{0}{x\cdot y}
&=&\matr{0}{x\cdot y}{0}{0}+\matr{0}{0}{0}{x\cdot y}\\
&=&2\tp{\matr{0}{y}{0}{0}}{\matr{v}{0}{0}{0}}{\matr{x}{0}{0}{0}}\\
&+&2\tp{\matr{0}{0}{0}{y}}{\matr{0}{0}{v}{0}}{\matr{0}{0}{x}{0}}
\end{eqnarray*}
From the last two displayed equations, we have
\[
\tp{\matr{0}{y}{0}{0}}{\matr{v}{0}{0}{0}}{\matr{x}{0}{0}{0}}
+\tp{\matr{0}{0}{0}{y}}{\matr{0}{0}{v}{0}}{\matr{0}{0}{x}{0}}\]
\[
=\tp{\matr{0}{0}{x}{0}}{\matr{v}{0}{0}{0}}{\matr{0}{y}{0}{0}}
+\tp{\matr{x}{0}{0}{0}}{\matr{0}{0}{v}{0}}{\matr{0}{0}{0}{y}}
\]
The first term on the left of the last displayed equation is of the form 
$\matr{0}{A}{0}{0}$ and the second is of the form 
$\matr{0}{0}{0}{B}$.
Again multiplying rows and columns by $-1$ and using the fact that isometries preserve the partial triple product shows that the first term on the right of the last displayed equation is of the form 
$\matr{0}{0}{0}{C}$ and the second is of the form 
$\matr{0}{D}{0}{0}$.
Thus 
\[
\tp{\matr{0}{y}{0}{0}}{\matr{v}{0}{0}{0}}{\matr{x}{0}{0}{0}}=
\tp{\matr{x}{0}{0}{0}}{\matr{0}{0}{v}{0}}{\matr{0}{0}{0}{y}}.\qed
\]

 \begin{lemma}\label{lem:1217111}
$x\cdot v=v\cdot x=x$ for every $x\in A$. 
 \end{lemma}
\pf\
Apply the main identity to write
\[
\tp
{\matr{0}{v}{0}{0}}
{\matr{0}{v}{0}{0}}
{\tp
{         \matr{v}{0}{0}{0}      }
{\matr{0}{v}{0}{0}}
{\matr{0}{x}{0}{0}}
}=R-S+T
\]
where
\begin{eqnarray*}
R&=&
\tp
{\tp
{\matr{0}{v}{0}{0}}
{\matr{0}{v}{0}{0}}
{\matr{v}{0}{0}{0}}
}
{{\matr{0}{v}{0}{0}}}
{{\matr{0}{x}{0}{0}}}\\
&=&\frac{1}{2}\tp
{\matr{v}{0}{0}{0}}
{{\matr{0}{v}{0}{0}}}
{{\matr{0}{x}{0}{0}}}
\end{eqnarray*}

\begin{eqnarray*}
S&=&
\tp
{\matr{v}{0}{0}{0}   }
{\tp
{\matr{0}{v}{0}{0}}
{\matr{0}{v}{0}{0}}
{\matr{0}{v}{0}{0}}
}
{{\matr{0}{x}{0}{0}}}\\
&=&\tp
{\matr{v}{0}{0}{0}}
{{\matr{0}{v}{0}{0}}}
{{\matr{0}{x}{0}{0}}}
\end{eqnarray*}
and
\begin{eqnarray*}
T&=&
\tp
{\matr{v}{0}{0}{0}   }
{\matr{0}{v}{0}{0}   }
{\tp
{\matr{0}{v}{0}{0}}
{\matr{0}{v}{0}{0}}
{\matr{0}{x}{0}{0}}
}\\
&=&\tp
{\matr{v}{0}{0}{0}}
{{\matr{0}{v}{0}{0}}}
{{\matr{0}{x}{0}{0}}}.
\end{eqnarray*}
Thus
\[
\tp
{\matr{0}{v}{0}{0}}
{\matr{0}{v}{0}{0}}
{\tp
{         \matr{v}{0}{0}{0}      }
{\matr{0}{v}{0}{0}}
{\matr{0}{x}{0}{0}}}=
\]
\begin{equation}\label{eq:1217111}
\frac{1}{2}\tp{\matr{v}{0}{0}{0}}
{\matr{0}{v}{0}{0}}
{\matr{0}{x}{0}{0}}
=\frac{1}{4}\matr{x\cdot v}{0}{0}{0}.
\end{equation}

Apply the main identity again to write
\[
\tp
{\matr{0}{v}{0}{0}}
{\matr{0}{v}{0}{0}}
{\tp
{         \matr{0}{v}{0}{0}      }
{\matr{0}{0}{0}{v}}
{\matr{0}{0}{x}{0}}
}=R'-S'+T'
\]
where
\[
R'=
\tp
{\tp
{\matr{0}{v}{0}{0}}
{\matr{0}{v}{0}{0}}
{\matr{0}{v}{0}{0}}
}
{{\matr{0}{0}{0}{v}}}
{{\matr{0}{0}{x}{0}}}
=\tp
{\matr{0}{v}{0}{0}}
{\matr{0}{0}{0}{v}}
{\matr{0}{0}{x}{0}},
\]
\[
S'=
\tp
{\matr{0}{v}{0}{0}   }
{\tp
{\matr{0}{v}{0}{0}}
{\matr{0}{v}{0}{0}}
{\matr{0}{0}{0}{v}}
}
{{\matr{0}{0}{x}{0}}}=
\frac{1}{2}\tp
{\matr{0}{v}{0}{0}}
{\matr{0}{0}{0}{v}}
{\matr{0}{0}{x}{0}},
\]
and
\[
T'=
\tp
{\matr{0}{v}{0}{0}   }
{\matr{0}{0}{0}{v}   }
{\tp
{\matr{0}{v}{0}{0}}
{\matr{0}{v}{0}{0}}
{\matr{0}{0}{x}{0}}
}=0
\]
by Lemma~\ref{lem:2.7}.
Thus
\[
\tp
{\matr{0}{v}{0}{0}}
{\matr{0}{v}{0}{0}}
{\tp
{         \matr{0}{v}{0}{0}      }
{\matr{0}{0}{0}{v}}
{\matr{0}{0}{x}{0}}}=R'-S'+T'=
\]
\begin{equation}\label{eq:1217112}
\frac{1}{2}\tp
{\matr{0}{v}{0}{0}}
{\matr{0}{0}{0}{v}}
{\matr{0}{0}{x}{0}}
=\frac{1}{4}\matr{v\cdot x}{0}{0}{0},
\end{equation}
the last step by Lemma~\ref{lem:2.18}.

By Lemmas~\ref{lem:2.8},\ref{lem:2.13} and \ref{lem:2.18}
\begin{equation}\label{eq:1217113}
\matr{x}{0}{0}{0}=\tp
{\matr{v}{0}{0}{0}}
{\matr{0}{v}{0}{0}}
{\matr{0}{x}{0}{0}}+
\tp
{\matr{0}{v}{0}{0}}
{\matr{0}{v}{0}{0}}
{\matr{0}{0}{x}{0}}.
\end{equation}

Adding (\ref{eq:1217111}) and (\ref{eq:1217112}) and using (\ref{eq:1217113}) results in 
\[
\frac{1}{2}\matr{v\cdot x}{0}{0}{0}=
\tp
{\matr{0}{v}{0}{0}}
{\matr{0}{v}{0}{0}}
{\matr{x}{0}{0}{0}}
=
\frac{1}{4}\matr{x\cdot v}{0}{0}{0}+\frac{1}{4}\matr{v\cdot x}{0}{0}{0}.
\]
Thus $v\cdot x=x\cdot v$ and since $x\cdot v+v\cdot x=2\tp{v}{v}{x}=2x$, the lemma is proved.\qed
\smallskip

We are now in a position to fill the gap in Lemma~\ref{lem:2.7}, which included only the case $x=v$ of Corollary~\ref{cor:4.4} below. To prove Corollary~\ref{cor:4.4}, we first need yet another lemma, along the lines of 
Lemmas~\ref{lem:2.5},\ref{lem:2.5prime},\ref{lem:2.10}, and \ref{lem:2.18}. We omit the by now standard proof, except to point out that the isometry involved is 
\[
\matr{0}{a}{0}{b}\mapsto \matr{0}{0}{a}{b}.
\]

\begin{lemma}\label{lem:1218111}
If $B,D\in A$ are defined by  $$\tp
{\matr{0}{a}{0}{b}}
{\matr{0}{v}{0}{0}}
{\matr{0}{c}{0}{d}}
=
\matr{0}{B}{0}{D},
$$
then
$$\tp
{\matr{0}{0}{a}{b}}
{\matr{0}{0}{v}{0}}
{\matr{0}{0}{c}{d}}
=
\matr{0}{0}{B}{D},
$$
In particular, 
\[
\tp
{\matr{0}{0}{0}{v}}
{\matr{0}{0}{v}{0}}
{\matr{0}{0}{x}{0}}
=
\tp
{\matr{0}{0}{0}{v}}
{\matr{0}{v}{0}{0}}
{\matr{0}{x}{0}{0}}
\]
\end{lemma}
\smallskip

\begin{corollary}\label{cor:4.4}
$\tp{\matr{y}{0}{0}{0}}{\matr{0}{0}{0}{v}}{\matr{0}{0}{0}{x}}=0.$
\end{corollary}
\pf\
By Lemma~\ref{lem:2.7} and the main identity,
\[
0=\tp
{\matr{0}{0}{0}{v}}
{\matr{0}{v}{0}{0}}
{\tp
{\matr{y}{0}{0}{0}}
{\matr{0}{0}{0}{v}}
{\matr{0}{x}{0}{0}}
}
=R-S+T
\]
where by two applications of  Lemma~\ref{lem:2.18}
\begin{eqnarray*}
R&=&
\tp
{\tp
{\matr{0}{0}{0}{v}}
{\matr{0}{v}{0}{0}}
{\matr{y}{0}{0}{0}}
}
{{\matr{0}{0}{0}{v}}}
{{\matr{0}{x}{0}{0}}}\\
&=&
\tp
{\tp
{\matr{0}{0}{y}{0}}
{\matr{0}{0}{0}{v}}
{\matr{0}{0}{0}{v}}
}
{{\matr{0}{0}{0}{v}}}
{{\matr{0}{x}{0}{0}}}\\
&=&\frac{1}{2}\tp
{\matr{0}{0}{v\cdot y}{0}}
{{\matr{0}{0}{0}{v}}}
{{\matr{0}{x}{0}{0}}}
=\frac{1}{2}\tp
{\matr{v\cdot y}{0}{0}{0}}
{{\matr{0}{v}{0}{0}}}
{{\matr{0}{x}{0}{0}}}\\
&=&\frac{1}{2}\tp
{\matr{y}{0}{0}{0}}
{{\matr{0}{v}{0}{0}}}
{{\matr{0}{x}{0}{0}}}
=\frac{1}{4}\matr{x\cdot y}{0}{0}{0},
\end{eqnarray*}
and by direct calculation
\begin{eqnarray*}
S&=&
\tp
{\matr{y}{0}{0}{0}   }
{\tp
{\matr{0}{v}{0}{0}}
{\matr{0}{0}{0}{v}}
{\matr{0}{0}{0}{v}}
}
{{\matr{0}{x}{0}{0}}}\\
&=&
\frac{1}{2}\tp
{\matr{y}{0}{0}{0}}
{\matr{0}{v}{0}{0}}
{\matr{0}{x}{0}{0}}
=\frac{1}{4}\matr{x\cdot y}{0}{0}{0}.
\end{eqnarray*}
Thus
\[
T=
\tp
{\matr{y}{0}{0}{0}   }
{\matr{0}{0}{0}{v}   }
{\tp
{\matr{0}{0}{0}{v}}
{\matr{0}{v}{0}{0}}
{\matr{0}{x}{0}{0}}
}=0.
\]
By Lemmas~\ref{lem:1217111} and ~\ref{lem:1218111},
\begin{eqnarray*}
0\ =\ T&=&
\tp
{\matr{y}{0}{0}{0}   }
{\matr{0}{0}{0}{v}   }
{\tp
{\matr{0}{0}{0}{v}}
{\matr{0}{0}{v}{0}}
{\matr{0}{0}{x}{0}}}
\\
&=&
\frac{1}{2}\tp
{\matr{y}{0}{0}{0}}
{\matr{0}{0}{0}{v}}
{\matr{0}{0}{0}{x\cdot v}}\\
&=&
\frac{1}{2}\tp
{\matr{y}{0}{0}{0}}
{\matr{0}{0}{0}{v}}
{\matr{0}{0}{0}{x}}. \qed
\end{eqnarray*}

The following proposition is critical.

\begin{proposition}\label{prop:2.14}
For $X,Y\in M_n(A)$, and $V=\hbox{ diag}\, (v,v,\ldots,v)=v\otimes I_n$,
\begin{description}
\item[(a)]
$
\{XVV\}=X
$
\medskip
\item[(b)]
$\left[\begin{array}{cc}
0& Y\cdot X\\
0&0
\end{array}
\right]=2
\left[\begin{array}{cc}
Y& 0\\
0&0
\end{array}
\right]
\left[\begin{array}{cc}
v\otimes I_n& 0\\
0&0
\end{array}
\right]
\left[\begin{array}{cc}
0& X\\
0&0
\end{array}
\right].
$
\medskip
\item[(c)]
$
X\cdot Y+Y\cdot X=2\{XVY\}.
$
\end{description}
\end{proposition}
\pf\
We shall prove by induction on $k$ that the proposition holds for $n=1,2,\ldots, 2k$.

If $n=1$,  (a) is the first assumption in Theorem~\ref{thm:main2}, (b) is Definition~\ref{def:3.13}, and (c) has been noted in (\ref{eq:0606121}) as a consequence of Lemmas~\ref{lem:2.8} and \ref{lem:2.14}.  \smallskip

Now let $n=2$.\footnote{Although the proof of this case is long, it renders the inductive step trivial}
Let us write
\begin{eqnarray*}
\lefteqn{\tp{\matr{v}{0}{0}{v}}{\matr{v}{0}{0}{v}}{\matr{a}{b}{c}{d}}=}\\
&&\tp{\matr{v}{0}{0}{0}}{\matr{v}{0}{0}{0}}{\matr{a}{b}{c}{d}}
+\tp{\matr{v}{0}{0}{0}}{\matr{0}{0}{0}{v}}{\matr{a}{b}{c}{d}}\\
&+&\tp{\matr{0}{0}{0}{v}}{\matr{v}{0}{0}{0}}{\matr{a}{b}{c}{d}}
+\tp{\matr{0}{0}{0}{v}}{\matr{0}{0}{0}{v}}{\matr{a}{b}{c}{d}}.
\end{eqnarray*}
The two middle terms on the right side of this equation vanish 
by Lemma~\ref{lem:2.7} and Corollary~\ref{cor:4.4}. 
The first term can be written as
\begin{eqnarray*}
\lefteqn{\tp{\matr{v}{0}{0}{0}}{\matr{v}{0}{0}{0}}{\matr{a}{b}{c}{d}}=}\\
&&\tp{\matr{v}{0}{0}{0}}{\matr{v}{0}{0}{0}}{\matr{a}{0}{0}{0}}
+\tp{\matr{v}{0}{0}{0}}{\matr{v}{0}{0}{0}}{\matr{0}{b}{0}{0}}\\
&+&\tp{\matr{v}{0}{0}{0}}{\matr{v}{0}{0}{0}}{\matr{0}{0}{c}{0}}
+\tp{\matr{v}{0}{0}{0}}{\matr{v}{0}{0}{0}}{\matr{0}{0}{0}{d}}\\
&=&\matr{\tp{v}{v}{a}}{0}{0}{0}
+\frac{1}{2}\matr{0}{v\cdot b}{0}{0}
+\frac{1}{2}\matr{0}{0}{v\cdot c}{0}
+\matr{0}{0}{0}{0}\\
&=&\matr{a}{b/2}{c/2}{0}.
\end{eqnarray*}
The last term can be written  (using Lemma~\ref{lem:1218111} in the second term) as
\begin{eqnarray*}
\lefteqn{\tp{\matr{0}{0}{0}{v}}{\matr{0}{0}{0}{v}}{\matr{a}{b}{c}{d}}=}\\
&&\tp{\matr{0}{0}{0}{v}}{\matr{0}{0}{0}{v}}{\matr{a}{0}{0}{0}}
+\tp{\matr{0}{0}{0}{v}}{\matr{0}{0}{0}{v}}{\matr{0}{b}{0}{0}}\\
&+&\tp{\matr{0}{0}{0}{v}}{\matr{0}{0}{0}{v}}{\matr{0}{0}{c}{0}}
+\tp{\matr{0}{0}{0}{v}}{\matr{0}{0}{0}{v}}{\matr{0}{0}{0}{d}}\\
&=&\matr{0}{0}{0}{0}
+\frac{1}{2}\matr{0}{v\cdot b}{0}{0}
+\frac{1}{2}\matr{0}{0}{v\cdot c}{0}
+\matr{0}{0}{0}{d}\\
&=&\matr{0}{b/2}{c/2}{d}.
\end{eqnarray*}

This completes the proof of (a) for $n=2$.
Lemmas~\ref{lem:2.5} to ~\ref{lem:2.14} and ~\ref{lem:2.18} to ~\ref{lem:1218111} now follow automatically for elements of $M_2(A)$,
 since the proofs for $M_2(A)$ are the same as those for $A$ once you have (a).\smallskip

Once (b) is proved for $n=2$, (c) will follow in the same way as (\ref{eq:0606121}) from the fact that Lemmas~\ref{lem:2.8} and \ref{lem:2.14}
are valid for $M_2(A)$.\smallskip

 The left side of (b) expands into 8 terms:
\begin{eqnarray*}
\left[\begin{array}{cc}
0& Y\cdot X\\
0&0
\end{array}
\right]&=&\matr{0}   {\matr    {y_{11}\cdot x_{11}}   {0}  {0}   {0}}{0}{0}
+\matr{0}{\matr{y_{12}\cdot x_{21}}{0}{0}{0}}{0}{0}\\
&+&
\matr{0}{    \matr{0}{y_{11}\cdot x_{12}}{0}{0}}{0}{0}
+\matr{0} {\matr{0}{y_{12}\cdot x_{22}}{0}{0}} {0} {0} \\
&+&
\matr{0}     {   \matr{0}{0}{ y_{21}\cdot x_{11}}{0}} {0} {0}
+\matr  {0}    {\matr{0}{0}{y_{22}\cdot x_{21}}{0}}    {0}     {0}\\
&+&
\matr{0}     {   \matr{0}{0}{0}{ y_{21}\cdot x_{12}}} {0} {0}
+\matr  {0}    {\matr{0}{0}{0}{y_{22}\cdot x_{22}}}    {0}     {0}
\end{eqnarray*}
For the right side, we have
\[
\tp
{\matr {\matr{y_{11}}{y_{12}}{y_{21}}{y_{22}}}{0}{0}{0}}
{\matr{\matr{v}{0}{0}{v}}{0}{0}{0}}
{\matr{0}{\matr{x_{11}}{x_{12}}{x_{21}}{x_{22}}}{0}{0}}
\]
which is the sum of 32 terms.  We show now that 24 of these 32 terms are zero, and each of 
the other 8 terms is equal to one of the 8 terms in the expansion of the left side.
We note first that by changing the signs of the first two columns we have that
\[
\tp
{\matr {\matr{y_{11}}{y_{12}}{y_{21}}{y_{22}}}{0}{0}{0}}
{\matr{\matr{v}{0}{0}{v}}{0}{0}{0}}
{\matr{0}{\matr{x_{11}}{x_{12}}{x_{21}}{x_{22}}}{0}{0}}
\]
has the form 
\[
\matr{0}{\matr{A}{B}{C}{D}}{0}{0}.
\]

We shall consider eight cases.\smallskip

Case 1A:  $Y=y_{11}\otimes e_{11}=\matr{y_{11}}{0}{0}{0},\ V=v\otimes e_{11}=\matr{v}{0}{0}{0}$\smallskip

In this case, further analysis shows that 
\[
\tp{\matr{\matr{y_{11}}{0}{0}{0}}{0}{0}{0} } {   \matr{\matr{v}{0}{0}{0}}{0}{0}{0}   } {    \matr{0}{\matr{x_{11}}{0}{0}{0}}{0}{0}   }\]
is of the form 
\[
\matr{0}{\matr{A}{0}{0}{0}}{0}{0}.
\]
and hence is unchanged by applying the isometry $C_{14}$ which interchanges the first and fourth columns.  The resulting (form of the) triple product we started with is therefore
\[
\tp {\matr{0}{\matr{0}{y_{11}}{0}{0}}    {0}{0} }    
      {\matr{0}{\matr{0}{v}{0}{0}}{0}{0}}
      {\matr{0}{\matr{x_{11}}{0}{0}{0}}{0}{0}}
      \]
      which equals (isometries preserve the partial triple product)
      \[
      \matr{0}{\tp{\matr{0}{y_{11}}{0}{0}}{\matr{0}{v}{0}{0}}{\matr{x_{11}}{0}{0}{0}}   }{0}{0}\\
=\frac{1}{2}\matr{0}{\matr{y_{11}\cdot x_{11}}{0}{0}{0}}{0}{0}
\]
as required.
An identical argument, using $C_{13}$ instead of $C_{14}$ shows that
\[
\tp{\matr{\matr{y_{11}}{0}{0}{0}}{0}{0}{0} } {   \matr{\matr{v}{0}{0}{0}}{0}{0}{0}   } {    \matr{0}{\matr{0}{x_{12}}{0}{0}}{0}{0}   }=\frac{1}{2}\matr{0}{\matr{0}{y_{11}\cdot x_{12}}{0}{0}}{0}{0}.
\]
To finish case 1A, use the isometry $R_{14}$ which interchanges the first and fourth rows on
\[
\tp
 { \matr{    \matr{y_{11}} {0} {0} {0} }   {0} {0} {0}      }
 { \matr{    \matr{v}{0}{0}{0}} {0}{0}{0}   }
 {\matr{0} {\matr{0}{0}{x_{21}}{x_{22}}}   {0}{0}   }
\]
to obtain
\[
\tp
 {       \matr{ 0} {0}{\matr{0}{0}{y_{11}}{0} } {0}    }
 {   \matr{ 0} {0}{\matr{0}{0}{v}{0} } {0}   }
 {\matr{0} {\matr{0}{0}{x_{21}}{x_{22}}}   {0}{0}   }
\]
which is zero by Lemma~\ref{lem:2.7}, which is valid for  $M_2(A)$.   Hence, the original triple product is zero.\medskip

Case 1B:  $Y=y_{11}\otimes e_{11}=\matr{y_{11}}{0}{0}{0},\ V=v\otimes e_{22}=\matr{0}{0}{0}{v}$\smallskip

Using the isometry $C_2(i)$ of multiplication of the second column by the imaginary unit $i$
we have that
\[
\tp
{\matr{\matr{y_{11}}{0}{0}{0}}{0}{0}{0}}
 {    \matr       {\matr      {0} {0} {0} {v}    } {0} {0} {0}          } 
 {     \matr      {0}     {\matr     {x_{11}} {x_{12}} {x_{21}} {x_{22} }  }   {0}   {0}         }
\]
which is of the form 
\[
\matr{0}{\matr{A}{B}{C}{D}}{0}{0},
\]
is equal to a non-zero multiple of  its negative, and is thus zero.\medskip

Case 2A:  $Y=y_{12}\otimes e_{12}=\matr{0}{y_{12}}{0}{0},\ V=v\otimes e_{11}=\matr{v}{0}{0}{0}$\smallskip

Using the isometry $C_1(i)$ of multiplication of the first column by the imaginary unit $i$
we have that
\[
\tp
{\matr{\matr{0}{y_{12}}{0}{0}}{0}{0}{0}}
 {    \matr       {\matr      {v} {0} {0} {0}    } {0} {0} {0}          } 
 {     \matr      {0}     {\matr     {x_{11}} {x_{12}} {x_{21}} {x_{22} }  }   {0}   {0}         }
\]
which is of the form 
\[
\matr{0}{\matr{A}{B}{C}{D}}{0}{0},
\]
is equal to a non-zero multiple of  its negative, and is thus zero.\medskip

Case 2B:  $Y=y_{12}\otimes e_{12}=\matr{0}{y_{12}}{0}{0},\ V=v\otimes e_{22}=\matr{0}{0}{0}{v}$\smallskip

Using the isometry $R_{23}$  which interchanges rows 2 and 3
we have that
\[
\tp
{\matr{\matr{0}{y_{12}}{0}{0}}{0}{0}{0}}
 {    \matr       {\matr      {0} {0} {0} {v}    } {0} {0} {0}          } 
 {     \matr      {0}     {\matr     {x_{11}} {x_{12}} {0} {0 }  }   {0}   {0}         }=0 
\]
by Lemma~\ref{lem:2.7}, which is valid for  $M_2(A)$. 

  Using the isometry $C_{24}$ and Lemma~\ref{lem:2.18}, we have that
\[
\tp
{\matr{\matr{0}{y_{12}}{0}{0}}{0}{0}{0}}
 {    \matr       {\matr      {0} {0} {0} {v}    } {0} {0} {0}          } 
 {     \matr      {0}     {\matr     {0} {0} {x_{21}} {0 }  }   {0}   {0}         }=
\frac{1}{2} \matr{0}{\matr{y_{12}\cdot  x_{21}}{0}{0}{0}}{0}{0}
\]

Using the isometry $C_{23}$ and Lemma~\ref{lem:2.18}, we have that
\[
\tp
{\matr{\matr{0}{y_{12}}{0}{0}}{0}{0}{0}}
 {    \matr       {\matr      {0} {0} {0} {v}    } {0} {0} {0}          } 
 {     \matr      {0}     {\matr     {0} {0} {0} {x_{22} }  }   {0}   {0}         }=
\frac{1}{2} \matr{0}{\matr{0}{y_{12}\cdot  x_{22}}{0}{0}}{0}{0}.
\]
\medskip

Case 3A:  $Y=y_{21}\otimes e_{21}=\matr{0}{0}{y_{21}}{0},\ V=v\otimes e_{11}=\matr{v}{0}{0}{0}$\smallskip

Using the isometry $R_{13}$  which interchanges rows 1 and 3
we have that
\[
\tp
{\matr{\matr{0}{0}{y_{21}}{0}}{0}{0}{0}}
 {    \matr       {\matr      {v} {0} {0} {0}    } {0} {0} {0}          } 
 {     \matr      {0}     {\matr     {0}{0}{x_{21}} {x_{22}}  }   {0}   {0}         }=0 
\]
by Lemma~\ref{lem:2.7}, which is valid for  $M_2(A)$.

  Using the isometry $C_{14}$ and Lemma~\ref{lem:2.18}, we have that
\[
\tp
{\matr{\matr{0}{0}{y_{21}}{0}}{0}{0}{0}}
 {    \matr       {\matr      {v} {0} {0} {0}    } {0} {0} {0}          } 
 {     \matr          {\matr      {x_{11}}{0}{0} {0 }}{0}     {0}   {0}         }=
\frac{1}{2} \matr{0}{\matr{0}{0}{y_{21}\cdot  x_{11}}{0}}{0}{0}
\]

Using the isometry $C_{13}$ and Lemma~\ref{lem:2.18}, we have that
\[
\tp
{\matr{\matr{0}{0}{y_{21}}{0}}{0}{0}{0}}
 {    \matr       {\matr      {v} {0} {0} {0}    } {0} {0} {0}          } 
 {     \matr      {0}     {\matr     {0}{x_{12}}  {0} {0 }  }   {0}   {0}         }=
\frac{1}{2} \matr{0} {       \matr{0}{0}{0}{y_{21}\cdot  x_{11}   }   }     {0}       {0}.
\]
\medskip

Case 3B:  $Y=y_{21}\otimes e_{21}=\matr{0}{0}{y_{21}}{0},\ V=v\otimes e_{22}=\matr{0}{0}{0}{v}$\smallskip

Using the isometry $C_2(i)$ of multiplication of the second column by the imaginary unit $i$
we have that
\[
\tp
{\matr{\matr{0}{0}{y_{21}}{0}}{0}{0}{0}}
 {    \matr       {\matr      {0} {0} {0} {v}    } {0} {0} {0}          } 
 {     \matr      {0}     {\matr     {x_{11}} {x_{12}} {x_{21}} {x_{22} }  }   {0}   {0}         }
\]
which is of the form 
\[
\matr{0}{\matr{A}{B}{C}{D}}{0}{0},
\]
is equal to a non-zero multiple of  its negative, and is thus zero.\medskip

Case 4A:  $Y=y_{22}\otimes e_{22}=\matr{0}{0}{0}{y_{22}},\ V=v\otimes e_{11}=\matr{v}{0}{0}{0}$\smallskip

Using the isometry $C_1(i)$ of multiplication of the first column by the imaginary unit $i$
we have that
\[
\tp
{\matr{\matr{0}{0}{0}{y_{22}}}{0}{0}{0}}
 {    \matr       {\matr      {v} {0} {0} {0}    } {0} {0} {0}          } 
 {     \matr      {0}     {\matr     {x_{11}} {x_{12}} {x_{21}} {x_{22} }  }   {0}   {0}         }
\]
which is of the form 
\[
\matr{0}{\matr{A}{B}{C}{D}}{0}{0},
\]
is equal to a non-zero multiple of  its negative, and is thus zero.\medskip

Case 4B:  $Y=y_{22}\otimes e_{22}=\matr{0}{0}{0}{y_{22}},\ V=v\otimes e_{22}=\matr{0}{0}{0}{v}$\smallskip

Using the isometry $R_{23}$ shows that
\[
\tp{\matr{\matr{0}{0}{0}{y_{22}}}{0}{0}{0} } {   \matr{\matr{0}{0}{0}{v}}{0}{0}{0}   } {    \matr{0}{\matr{x_{11}}{x_{12}}{0}{0}}{0}{0}   }=0.
\]

Using the isometry $C_{24}$ shows that
\[
\tp
{\matr{\matr{0}{0}{0}{y_{22}}}{0}{0}{0}}
 {    \matr       {\matr      {0} {0} {0} {v}    } {0} {0} {0}          } 
 {     \matr          {\matr      {0}{0}{x_{21}} {0 }  }   {0}   {0}   {0}      }=
\frac{1}{2} \matr{0}{\matr{0}{0}{y_{22}\cdot  x_{21}}{0}}{0}{0}
\]

Using the isometry $C_{23}$ shows that
\[
\tp
{\matr{\matr{0}{0}{0}{y_{22}}}{0}{0}{0}}
 {    \matr       {\matr      {0} {0} {0} {v}    } {0} {0} {0}          } 
 {     \matr          {\matr      {0}{0}{0} {x_{22} }  }   {0}   {0}    {0}     }=
\frac{1}{2} \matr{0}{\matr{0}{0}{0}{y_{22}\cdot  x_{22}}}{0}{0}
\]

This completes the proof of (b) and hence of (c) for $n=2$, and the proposition for $k=1$.
\smallskip

We now assume the the proposition holds for $ n=1,2,\ldots,2k$.  For  any $X\in M_n(A)$, let us write
\[
\tilde X=\left\{\begin{array}{cl}
\matr{X}{0}{0}{0}& \hbox{ if }n=2k+1,\\
\medskip
X& \hbox{ if }n=2k+2.
\end{array}\right.
\]
We then write
\[
\tilde X=\matr{X_{11}}{X_{12}}{X_{21}}{X_{22}}
\]
where $X_{ij}\in M_{k+1}(A)$.
 Since $k+1\le 2k$, the induction proceeds by simply repeating  the proofs in the case $n=2$, with $X,Y,V$ replaced by $\tilde X,\tilde Y, \tilde V$.
\qed
\smallskip

We can now complete the proof of  our main result.

\begin{theorem}\label{thm:main2}
An operator space $A$ is completely isometric to a unital operator
algebra if and only there exists $v\in CS(A)$ such that:
\begin{description}
\item[(i)] $h_v(x+v)-h_v(x)-h_v(v)+v=-2x$ for all $x\in A$
\item[(ii)] Let  $V=\hbox{diag}(v,\ldots,v)\in M_n(A)$.  For all $X\in M_n(A)$
\[
\|V-h_{V}(X)\|\le \|X\|^2.
\]
\end{description}
\end{theorem}
\pf\
As was already pointed out, the first assumption is equivalent to the condition $\tp{x}{v}{v}=x$, so that all the machinery developed so far is available.  In particular, $v$ is a unit element for the product $x\cdot y$ and 
  for every $X\in M_2(A)$, $X\cdot X=\tp{X}{V}{X}$.

With $X=\matr{0}{x}{y}{0}$ for elements $x,y\in A$ of norm 1, we have 
\begin{eqnarray*}
\max(\|x\cdot y\|,\|y\cdot x\|)&=&\left\|\matr{x\cdot y}{0}{0}{y\cdot x}  \right\|=\|X\cdot X\|=\|\{XVX\|\|\\
&\le& \|X\|^2=\left\|\matr{0}{x}{y}{0}\right\|^2=\max(\|x\|,\|y\|)^2=1
\end{eqnarray*}
so the multiplication on $A$ is contractive.   The same argument shows that if $X,Y\in M_n(A)$, then $\|X\cdot Y\|\le \|X\|\|Y\|$ so the multiplication is completely contractive.  The result now follows from \cite{BleRauSin90}.

For the sake of completeness, we include the detail of the last inequality:
\begin{eqnarray*}
\max(\|X\cdot Y\|,\|Y\cdot X\|)&=&\left\|\matr{X\cdot Y}{0}{0}{Y\cdot X}  \right\|\\
&=&\left\|
\matr{0}{X}{Y}{0}\cdot \matr{0}{X}{Y}{0}\right\|\\
&=&\|\tp{\matr{0}{X}{Y}{0}}{\matr{V}{0}{0}{V}}{\matr{0}{X}{Y}{0}}\|\\ &\le&\|\matr{0}{X}{Y}{0}\|^2=\max(\|X\|,\|Y\|)^2.\qed
\end{eqnarray*}

\smallskip

\begin{remark}\label{rem:4.6}
The second condition in  Theorem~\ref{thm:main2} can be replaced by the following.\smallskip

(ii$^\prime$) Let $\tilde V$ denote the $2n$ by $2n$ matrix $\matr{V}{0}{0}{0}$, where $V=\hbox{diag}(v,\ldots,v)\in M_n(A)$.  For all $X,Y\in M_n(A)$
\[
\|h_{\tilde V}(\matr{Y}{X}{0}{0})-h_{\tilde V}(\matr{0}{X}{0}{0})-h_{\tilde V}(\matr{Y}{0}{0}{0}+\tilde V\|\le \|X\|\|Y\|.
\]
Equivalently\footnote{Although the 1/2 in (\ref{eq:0605121}) conveniently cancels the 2 in Proposition~\ref{prop:2.14}(b), its presence is justified by the fact that (\ref{eq:0605121}) holds in case $A$ is an operator algebra}, \begin{equation}\label{eq:0605121}
\|\tp{\matr{Y}{0}{0}{0}}{\matr{V}{0}{0}{0}}{\matr{0}{X}{0}{0}}\|\le\frac{1}{2}\|X\|\|Y\|,
\end{equation}
so by Proposition~\ref{prop:2.14}
\[
\|Y\cdot X\|\le \|X\|\|Y\|.
\]
By Lemma~\ref{lem:1217111} and the first condition,  $A$ is a unital (with a unit of norm 1 and not necessarily associative) algebra.   Remark~\ref{rem:4.6}  now follows from \cite{BleRauSin90}.
\end{remark}

\medskip

\medskip

We close by stating two problems for Banach spaces and three problems for operator spaces which arose in connection with this paper.

\begin{problem}\label{prob:1}
Is there a Banach space with partial triple product $\{x,a,y\}$ for which the inequality 
\[
\|\{x,a,y\}\|\le \|x\|\|a\|\|y\|
\]
does not hold?
\end{problem}

\begin{problem}\label{prob:2}
Is the symmetric part of the predual of a von Neumann algebra equal to 0?  What about the predual of a $JBW^*$-triple which does not contain  a Hilbert space as a direct summand?
\end{problem}

\begin{problem}\label{prob:3}
Is the completely symmetric part of an infinite dimensional Cartan factor of type 2,3 or 4 zero, as in the finite dimensional case?
\end{problem}

It is clear that the intersection of an operator space $A$ with the quasimultipliers of $I(A)$ from \cite{Kaneda04} (referenced in the second paragraph of 1.1 above) is a TRO and is contained in the completely symmetric part of $A$.  Our main theorem is that certain elements in the holomorphically defined completely symmetric part induce operator algebra products on $A$ while \cite{Kaneda04} shows that all operator algebra products on $A$  arise from the more concretely and algebraically defined quasimultpliers.  Hence it is natural to ask 

\begin{problem}\label{prob:5} Under what conditions does the completely symmetric part of an operator space  consist of quasimultipliers?
\end{problem}

Of course using direct sums and the discussion in the last two paragraphs of section 1.2, we can construct operator spaces whose completely symmetric part is different from zero and from the symmetric part of the operator space. However it would be more satisfying to answer the following problem.

\begin{problem}\label{prob:4}
Is there an operator space whose completely symmetric part is not contractively complemented, different from zero, and different from the symmetric part of the operator space?
\end{problem}

\bibliographystyle{amsplain}

\end{document}